\documentclass[]{sn-jnl}
\usepackage[utf8]{inputenc}
\usepackage{amsmath,amsthm,amsfonts,mathtools}
\usepackage{xspace}
\usepackage{enumitem}
\usepackage[numbers,sort&compress]{natbib}
\usepackage{tabularx}
\usepackage{arydshln}
\usepackage{array,multirow,graphicx,adjustbox}
\usepackage{xcolor}
\newcolumntype{L}[1]{>{\raggedright\arraybackslash}X{#1}}

\newcommand{\setN}{\mathcal{N}}
\newcommand{\setP}{\mathcal{P}}
\newcommand{\setM}{\mathcal{M}}
\newcommand{\setF}{\mathcal{F}}
\newcommand{\setR}{\mathcal{R}}
\newcommand{\setA}{\mathcal{A}}
\newcommand{\setS}{\mathcal{S}}
\newcommand{\setEarly}{\mathcal{S}^{\textnormal{early}}}
\newcommand{\setLate}{\mathcal{S}^{\textnormal{late}}}
\newcommand{\setNight}{\mathcal{S}^{\textnormal{night}}}

\newcommand{\setL}{\mathcal{L}}
\newcommand{\setE}{\mathcal{E}}
\newcommand{\setNprev}{\mathcal{N}^{\textnormal{prev}}}

\newcommand{\assign}{\textnormal{assign}}
\DeclareMathOperator{\skilllevel}{\textnormal{skill\_level}}
\DeclareMathOperator{\skillreq}{\textnormal{skill\_req}}
\DeclareMathOperator{\wload}{\textnormal{w\_load}}
\DeclareMathOperator{\maxload}{\textnormal{max\_load}}

\newcommand{\numbeds}{\textnormal{num\_beds}}
\DeclareMathOperator{\dist}{\textnormal{dist}}

\DeclareMathOperator{\adshift}{\textnormal{ad\_shift}}
\DeclareMathOperator{\dishift}{\textnormal{di\_shift}}

\DeclareMathOperator{\wpat}{\textnormal{walk\_pat}}
\newcommand{\age}{\textnormal{age}}
\DeclareMathOperator{\agegroup}{\textnormal{age\_group}}

\newcommand{\prev}{\textnormal{prev}}

\newcommand{\trans}{\textnormal{trans}}
\newcommand{\vio}{\textnormal{vio}}
\newcommand{\load}{\textnormal{load}}
\newcommand{\fair}{\textnormal{fair}}
\newcommand{\skill}{\textnormal{skill}}
\newcommand{\viogender}{\textnormal{gender}}
\newcommand{\maleinroom}{\textnormal{m\_in\_room}}
\newcommand{\femaleinroom}{\textnormal{f\_in\_room}}
\newcommand{\inroom}{\textnormal{in\_room}}
\newcommand{\bothrooms}{\textnormal{both\_rooms}}
\newcommand{\everassigned}{\textnormal{ever\_assigned}}

\newcommand{\skillnurses}{\textnormal{skill\_nurses}}
\newcommand{\maxshifts}{\textnormal{max\_shifts}}

\newcommand{\contribution}{\textnormal{ContribTable}}
\newcommand{\hetmatrix}{\textnormal{HetMatrix}}
\newcommand{\setNcomb}{\mathcal{N}^{\textnormal{comb}}}
\begin{document}

\title[Integrated patient-to-room and nurse-to-patient assignment]{Integrated patient-to-room and nurse-to-patient assignment in hospital wards}

\author[1]{\fnm{Tabea} \sur{Brandt}}\email{brandt@combi.rwth-aachen.de}

\author[2]{\fnm{Tom Lorenz} \sur{Klein}}\email{tom.klein@tum.de}

\author[3,4]{\fnm{Melanie} \sur{Reuter-Oppermann}}\email{m.n.reuter-oppermann@utwente.nl}

\author*[5]{\fnm{Fabian} \sur{Schäfer}}\email{fab.schaefer@tum.de}

\author[2,6]{\fnm{Clemens} \sur{Thielen}}\email{clemens.thielen@tum.de}

\author[3,7]{\fnm{Maartje} \sur{van de Vrugt}}\email{n.vandevrugt@amsterdamumc.nl}

\author[8]{\fnm{Joe} \sur{Viana}}\email{joe.viana@bi.no}

\affil[1]{\orgdiv{Combinatorial Optimization}, \orgname{RWTH Aachen University}, \orgaddress{\street{Ahornstr.~55}, \postcode{52074}~\city{Aachen}, \country{Germany}}}

\affil[2]{\orgdiv{TUM Campus Straubing for Biotechnology and Sustainability}, \orgname{Weihenstephan-Triesdorf University of Applied Sciences}, \orgaddress{\street{Am~Essigberg~3}, \postcode{94315}~\city{Straubing}, \country{Germany}}}

\affil[3]{\orgdiv{Center for Healthcare Operations Improvement \& Research}, \orgname{University of Twente}, \orgaddress{\street{Drienerlolaan 5}, \postcode{7522~NB}~\city{Enschede}, \country{The~Netherlands}}}

\affil[4]{\orgdiv{Information Systems, Software \& Digital Business Group}, \orgname{Technical University of Darmstadt}, \orgaddress{\street{Hochschulstr.~1}, \postcode{64289}~\city{Darmstadt}, \country{Germany}}}

\affil[5]{\orgdiv{Chair of Supply and Value Chain Management}, \orgname{Technical University of Munich}, \orgaddress{\street{Am~Essigberg~3}, \postcode{94315}~\city{Straubing}, \country{Germany}}}

\affil[6]{\orgdiv{Department of Mathematics, School of Computation, Information and Technology}, \orgname{Technical University of Munich}, \orgaddress{\street{Bolzmannstr.~3}, \postcode{85748}~\city{Garching bei München}, \country{Germany}}}

\affil[7]{\orgdiv{Department of Care Support, Strategy and Innovation}, \orgname{Amsterdam University Medical Centers}, \orgaddress{\street{De~Boelelaan~1117}, \postcode{1081~HV}~\city{Amsterdam}, \country{The~Netherlands}}}

\affil[8]{\orgdiv{Department of Accounting and Operations Management}, \orgname{BI~Norwegian Business School}, \orgaddress{\street{Nydalsveien~37}, \postcode{N-0484}~\city{Oslo}, \country{Norway}}}
\enlargethispage{0.65\baselineskip}

%%%%%ORCID of authors
%Joe Viana: 0000-0001-6018-8242
%Tabea Brandt: 0000-0002-8252-1891
% Clemens Thielen 0000-0003-0897-3571
% Maartje van de Vrugt 0000-0003-4287-1165
%Melanie Reuter-Oppermann 0000-0003-2231-7749

\abstract{
Assigning patients to rooms and nurses to patients are critical tasks within hospitals that directly affect patient and staff satisfaction, quality of care, and hospital efficiency. Both patient-to-room assignments and nurse-to-patient assignments are typically agreed upon at the ward level, and they interact in several ways such as jointly determining the walking distances nurses must cover between different patient rooms. This motivates to consider both problems jointly in an integrated fashion.

This paper presents the first optimization models and algorithms for the integrated patient-to-room and nurse-to-patient assignment problem. We provide a mixed integer programming formulation of the integrated problem that considers the typical objectives from the single problems as well as additional objectives that can only be properly evaluated when integrating both problems. Moreover, motivated by the inherent complexity that results from integrating these two \textsf{NP}-hard and already computationally challenging problems, we devise an efficient heuristic for the integrated patient-to-room and nurse-to-patient assignment problem. To evaluate the running time and quality of the solution obtained with the heuristic, we conduct extensive computational experiments on both artificial and real-world instances. The artificial instances are generated by a parameterized instance generator for the integrated problem that is made freely available.}

% Old version (Fabian's ORAHS abstract)
% The COVID-19 pandemic has highlighted the crucial role of healthcare staff, particularly in hospitals, where many countries have experienced staff shortages even before the pandemic and are still experiencing shortages. Therefore, efficient staff scheduling that considers staff preferences and satisfaction is a crucial task for hospitals. The number of beds typically measures hospital size, but the number of staffed beds is more relevant to patient care. Assigning patients to rooms and nurses to patients are critical tasks that directly affect staff satisfaction, quality of care, and hospital efficiency. Nurse-to-patient assignments are typically agreed upon at the ward level at the start of each shift, making it reasonable to plan them jointly with the patient-to-room assignments on the same ward. Therefore, we present an integrated optimization model for patient-to-room and nurse-to-patient assignments to enhance staff satisfaction and patient care. However, the resulting integrated planning problem is inherently complex since even the subproblems for patient-to-room and nurse-to-patient assignments are known to be NP-hard and computationally challenging. To evaluate the model's practical applicability, we conduct experiments to assess run time and solution quality for various settings, including different numbers of patients, shifts, and rooms.


\maketitle

\section{Introduction}\label{sec:introduction}
% Motivate why efficient planning of resources in hospitals is necessary.
For many years, an ever-rising demand for healthcare and increasing healthcare expenditures challenge hospitals to increase the efficiency of their operations~\cite{Drupsteen+etal:integrative-practices}. This results in a need for advanced managerial planning approaches that help to use the available scarce resources as efficiently as possible. Consequently, a wide range of methods and approaches have been developed in the Operations Research (OR) literature that aims at improving resource utilization through efficient planning~\cite{Rais+Viana:healt-care-or-survey,Hulshof:taxonomy,Jha:survey}. In particular, quantitative decision support has been proposed for important resources such as operating rooms~\cite{vanriet2015:trade-offs-or,guerriero2011:survey-or}, intensive care units~\cite{Bai+etal:intensive-care-review}, inpatient beds~\cite{He2019:sys-review-bed-management}, physicians~\cite{Erhard+etal:phys-sched-survey}, and nurses~\cite{Benazzouz+etal:survey,Clark+etal:nurse-rescheduling-survey}. 

\medskip

%Motivate integrated planning of several resources in hospitals
While efficient planning of single resources can already lead to improved resource utilization and large efficiency gains, it ignores the inherent complex interactions between different resources~\cite{Hulshof:taxonomy} and, as a consequence, might lead to suboptimal decisions on a system level. Therefore, a need for OR models and methods for \emph{integrated planning} of several resources has been identified~\cite{Hulshof:taxonomy,Jun+etal:discrete-event-survey,vanberkelSurvey}. This need is particularly apparent in hospitals, where different resources are typically required and used for treating patients. A recent literature review on integrated planning of multiple resources in hospitals is provided in~\cite{Rachuba+etal:hospital-review}.

\medskip

% Motivate beds and staff as critical resources in hospitals and introduce the PRA and NPA problems
Rooms and beds are critical assets of hospitals since they account for a considerable part of a hospital’s infrastructure and large financial investments are necessary for equipping them with medical devices that facilitate patient care~\cite{Vancroonenburg2016}. Additionally, a shift in demographics, the growing number of patient admissions, and rising inpatient units costs lead to high overall bed occupancy levels and require an increased focus on efficient bed management~\cite{Schaefer2019,He2019:sys-review-bed-management}. 
On the operational level, an important planning problem typically referred to as \emph{patient-to-room assignment~(PRA)} consists of assigning patients to suitable rooms such that a variety of constraints concerning, for instance, the patient's medical needs (e.g., required medical equipment) and preferences (e.g., concerning age and gender of roommates) are satisfied, while available room capacities are respected and transfers of patients between rooms are avoided~\cite{Demeester2010,Ceschia2011,Ceschia2012,Schaefer2019}. Different variants of the PRA problem have been studied extensively in the literature both in the static setting, when all information about patients and their admission and discharge times is known in advance, and in dynamic settings, when new patients may arrive unexpectedly -- see Section~\ref{subsec:literature-PRA} for a detailed overview.

\medskip

Medical staff also represent a particularly critical resource in hospitals since (1)~medical staff are involved in most patient-related activities in a hospital, and (2)~medical staff are a particularly scarce resource due to a general shortage of nurses~\cite{Aiken+etal:nurse-burnout} and physicians~\cite{Bodenheimer+Smith:primary-care,Erhard+etal:phys-sched-survey,Thielen2018}. This has led to increasing workloads for the staff over the last decades and makes good planning of medical staff a central concern for hospitals~\cite{Benazzouz+etal:survey}. In particular, aspects such as a fair distribution of workloads among staff members have a large impact on employee satisfaction and the efficient operation of a hospital. Concerning the nursing staff, distributing work fairly among nurses is considered essential for optimal quality of care~\cite{Mullinax2002}. 
Here, the workload of each single nurse is mostly determined by the patients the nurse is assigned to and their care requirements. %, but also by the walking distances that result from traveling between their assigned patients' rooms and other locations such as the nearest nursing station~\cite{Acar2016}. 
Consequently, determining a suitable \emph{nurse-to-patient assignment~(NPA)} that balances the workloads of the nurses represents an important operational problem that has received considerable attention in the literature -- see Section~\ref{subsec:literature-NPA} for an extensive literature review.

\medskip

% Motivate the integration of PRA and NPA - TODO!
With a few exceptions detailed in Section~\ref{subsec:literature-integration}, the PRA problem and the NPA problem have mostly been considered separately in the literature -- although there are important interactions between them. For instance, studies show that the walking distances that result from traveling between their assigned patients' rooms and other locations such as the nearest nursing station have a substantial impact on a nurse's workload during a shift~\cite{Acar2016,Butt2004}. These walking distances, however, can only be determined and optimized when considering PRAs and NPAs jointly. Moreover, it has been observed that assigning the minimum possible number of nurses to patients in the same room can help to minimize negative effects such as the spread of infections between rooms by nurses or the disturbance of patients by other nurses entering their room~\cite{Halwani2006,Cohen2012,Eveillard2009,Dancer2009}. This provides a strong motivation for integrating the PRA problem and the NPA problem by considering them jointly in one optimization model.

\medskip

% Detailed description of our contribution
While several publications motivate and discuss the integration of the two assignment problems (see Section~\ref{subsec:literature-integration}), this paper explicitly considers decisions on PRAs and NPAs in one integrated optimization problem for the first time. Besides the objectives classically considered in the two separate problems, this integrated problem also allows the evaluation of additional objectives that rely on the interaction of PRAs and NPAs. Based on existing studies on nurse workloads~\cite{Acar2016,Butt2004}, these objectives include the walking distances of nurses between assigned patients' rooms and additional relevant rooms such as the nearest nursing station. Moreover, also motivated by findings from the literature~\cite{Halwani2006,Cohen2012,Eveillard2009}, assigning the minimum possible number of nurses to patients in the same room is considered as an objective in order to mitigate negative effects such as the spread of infections between rooms by nurses or the disturbance of patients by other nurses entering their room.

% \medskip

In order to solve the computationally challenging integrated PRA and NPA problem, we provide a formulation of the problem as a mixed integer program (MIP) as well as an efficient heuristic. The heuristic extends the heuristic for the PRA problem presented in~\cite{Schaefer2019} to the integrated problem and additionally employs a new heterogeneity check between patient admission and discharge times for the room assignment part. Both the MIP and the heuristic are evaluated in extensive experimental results on real-world instances obtained from a ward of our partner hospital Amsterdam University Medical Centers (Amsterdam, The Netherlands) as well as artificial instances. The artificial instances are generated by a parameterized instance generator for the integrated problem that is made freely available. While the MIP only addresses the static version of the integrated problem in which all information about patients is known in advance, the heuristic can be easily adapted to dynamic settings where new patients arrive after some assignments have already been fixed.

\medskip

% Structure of the paper
\enlargethispage{\baselineskip}
The remainder of the paper is structured as follows. Section~\ref{sec:literature} summarizes the related literature on the PRA and the NPA problems as well as existing work related to their integration. Afterward, Section~\ref{sec:probDef} introduces the integrated PRA and NPA problem, while Section~\ref{sec:MIP} presents our MIP formulation of the problem. Afterward, Section~\ref{sec:solution-methods} presents a sequential solution approach based on a natural decomposition of the MIP formulation as well as our heuristic for the problem. 
Section~\ref{sec:instances} then describes both the developed instance generator for generating artificial test instances as well as the real-world instances obtained from a ward of our partner hospital. Our experimental results obtained on both types of instances are presented in Section~\ref{sec:exp_results}. The paper concludes in Section~\ref{sec:conclusion} with a summary and an outlook on future research.

\section{Related literature}\label{sec:literature}

% Patient-to-room assignment
% -> Tabea schreibt Entwurf

In this section, we first summarize the state of the art concerning PRA and NPA separately before discussing existing work on the integration of the two problems.

\subsection{Patient-to-room assignment}\label{subsec:literature-PRA}

The static version of the PRA problem has been formally introduced by Demeester et al.~\cite{Demeester2010} under the name of patient admission scheduling problem.\footnote{We refer to the problem exclusively as the PRA problem in the following since this is the most common term used in the recent literature. Moreover, the original term \emph{patient admission scheduling problem} is also used with a different meaning in the literature.}
In this version, all information about patients is known in advance, and the task is to assign patients to suitable rooms such that room capacity and gender policy are respected while minimizing both patient transfers and penalty costs for undesirable PRAs.
One important characteristic in this definition is that not every patient can be assigned to every room and patients may also have preferences towards specific rooms based on, e.g., available equipment or the number of beds.

\medskip

For the static PRA problem, mostly heuristic solution approaches are proposed in the literature, e.g.,
a tabu search algorithm~\cite{Demeester2010}, a local search algorithm~\cite{Ceschia2011}, a destroy and repair matheuristic~\cite{Guido2018}, and algorithms based on the Hungarian algorithm~\cite{Borchani2021}, column generation~\cite{Range2014}, or MIP~\cite{Thuran2017}.
Currently the best solutions for the benchmark instances provided by Demeester et al.~\cite{Demeester2010} are found by the heuristic of Guido et al.~\cite{Guido2018} and
by the exact, MIP-based solution approach proposed by Bastos et al.~\cite{Bastos2019}.
However, the exact approach uses considerably more computation time.

\medskip

%Further papers that consider this variant: Range~\cite{Range2014}, Thuran and Bilgen (2017)~\cite{Thuran2017}, Guido et al. (2018)~\cite{Guido2018} (best results for Demeester instances until 2019), Bastos et al. (2019)~\cite{Bastos2019}, Borchani et al. (2021)~\cite{Borchani2021} (use Demeester instances to test their methods)

The complexity of the static PRA problem is studied by Vancroonenburg et al.~\cite{Vancroonenburg2014} using its correspondence to the red-blue-transportation problem.
They show that the PRA problem is $\mathcal{NP}$-hard in general and even if all rooms have a capacity of three.
This result is strengthened by Ficker et al.~\cite{Ficker2021}, who prove that the PRA problem is also $\mathcal{NP}$-hard if all rooms have capacity two.

\medskip

%Taramasco et al. studied a similar but slightly different problem definition that includes the possibility of assigning a patient to a bed in a different hospital and propose an autonomous bat algorithm~\cite{Taramasco2019}.
%Maybe don't cite this? Seems like they consider only two days?

%Uncertainties and dynamic version
Ceschia and Schaerf~\cite{Ceschia2011, Ceschia2012} propose a dynamic version of the PRA problem that includes the handling of emergency patients and uncertainty in the patients' length-of-stay.
They propose a metaheuristic based on simulated annealing and an instance generator as well as a set of benchmark instances.
Vancroonenburg et al.~\cite{Vancroonenburg2016} studied a similar problem version using two online integer linear programming (ILP) models.
Lusby et al.~\cite{Lusby2016} propose a large neighborhood search heuristic for the dynamic PRA problem as proposed by Ceschia and Schaerf.

\medskip

%Uncertainties / dynamic version that introduces further new aspects:
A different approach for incorporating emergency patients is taken by
Schäfer et al.~\cite{Schaefer2017,Schaefer2019,Schaefer2021}, who use a rolling horizon approach, i.e., recomputation of the solution whenever a new event occurs, using a fast heuristic.
%and compares final solution to optimal solution that could have been obtained if all input data would have been available from the start. New aspects:
In their problem definition, they consider objectives for three stakeholders (patients, nurses, and physicians) and they are also the first to consider interdependencies between patients in the same room. % (fit between patients in the same room), and it is preferred if a physician is assigned to several patients in the same room (similar to what we do for the nurses).

\medskip

For a more detailed overview of the different versions and solution approaches for the PRA problem, we refer to Zhu et al.~\cite{Zhu2019}, who study the compatibility of short-term and long-term objectives in the context of dynamic PRA.

% \medskip

%%%%%%%%%%%%%%%%%%%%%%%%%%%%%%%%%%%%%%%%%%%%%%%%%%%
% Nurse-to-patient assignment
\subsection{Nurse-to-patient assignment}\label{subsec:literature-NPA}

The NPA problem is also considered by many different authors, where the most common objective is balancing the workloads of the nurses. For instance, Mullinax and Lawley~\cite{Mullinax2002} use this objective in the daily assignment of newborn infants to nurses in an intensive care nursery that is divided into several zones (rooms). Here, each infant might yield a different workload depending on their acuity and each nurse can be assigned a certain maximum number of infants to take care of at the beginning of the shift, but all of these need to be from the same zone. Since they find the problem to be too hard to solve using an integer program, they present a two-step heuristic approach that exploits the subdivision into zones by first computing the number of nurses allocated to each zone before assigning patients to nurses in each zone independently. 

\medskip

The problem introduced by Mullinax and Lawley~\cite{Mullinax2002} is also considered in~\cite{Sir2015,Ku2014,Schaus2014,Marzouk2020} -- each time with a slightly different objective function. Based on the formulation proposed in~\cite{Mullinax2002}, Sir et al.~\cite{Sir2015} formulate four MIP models for NPA that model the workloads of nurses by either the patient acuity indicators from a patient classification system (PCS), survey-based nurse-specific workload scores, or a combination of the two. Ku et al.~\cite{Ku2014} focus on minimizing the variance of the nurses' workloads using mixed integer quadratic programming (MIQP) and constraint programming (CP), while Schaus and Régin~\cite{Schaus2014} minimize the variance using a two-step decomposition approach that first computes the number of nurses allocated to each zone (which is done optimally by solving a resource allocation problem) before assigning patients to nurses in each zone independently using CP. Finally, Marzouk and Kamoun~\cite{Marzouk2020} formulate a binary integer program that assigns nurses to zones and individual patients with the objective of minimizing the total number of nurses used in a shift.

\medskip

Other work on NPA includes Punnakitikashem et al.~\cite{Punnakitikashem2008}, who present a stochastic integer programming model with the objective of minimizing excess workload for nurses and compare their approach to several other assignment policies (random assignment without considering workload, a simple heuristic, and solving the mean value problem using a deterministic integer program). Sundaramoorthi et al.~\cite{Sundaramoorthi2009} then use three of the assignment policies from~\cite{Punnakitikashem2008} as well as a clustered assignment policy to test their developed simulation model for evaluating NPAs. 

\medskip

While most of the literature on NPA mentioned above uses patient acuity as the main factor influencing nursing workloads, Acar and Butt~\cite{Acar2016} perform a detailed study in order to identify the activities that comprise a nurse's workload. They find that nurses spend a substantial part of their time traveling (walking) between locations, where travel between patient rooms and the nursing station is the most common type of travel. Here, according to Butt et al.~\cite{Butt2004}, the distance traveled by nurses is correlated to their assigned patient load and location, and key distances influencing the total travel distance of a nurse are the distances between assigned patients' rooms and (1) the nearest nursing station, (2) the nearest supply room, and (3) other assigned patients' rooms. Still, walking distances of nurses have not yet been considered explicitly as an objective in the literature on NPA since their minimization requires the simultaneous optimization of PRAs. This also holds for the objective of assigning the minimum possible number of nurses to patients in the same room, although it is known that assigning all patients in the same room to the same nurse or a small pool of nurses helps to avoid the transfer of hospital-acquired infections~\cite{Halwani2006,Cohen2012,Eveillard2009}, in particular \emph{Methicillin-resistant Staphylococcus aureus} (MRSA)~\cite{Dancer2009}.

\subsection{Integration of PRA and nurse-to-patient assignment}\label{subsec:literature-integration}

While we are not aware of any papers that explicitly integrate decisions concerning PRA and NPA in one optimization model, there still exists some literature considering the interplay between the two problems or between related problems. For instance, Thomas et al.~\cite{Thomas2013} present a mixed-integer goal programming model for the PRA problem that takes nurses into account via constraints on the required nurse-to-patient ratio in each unit of a hospital. Bilgin et al.~\cite{Bilgin2012} develop a general, high-level hyper-heuristic approach that can be used for both the PRA problem and the nurse rostering problem. Pesant~\cite{Pesant2016} addresses the integration of the nurse staffing problem (assigning an appropriate number of nurses to each unit within a ward given a nurse roster) and the NPA problem in a neonatal intensive care unit by solving CP models for the two problems consecutively, and Punnakitikashem et al.~\cite{Punnakitikashem2013} extend the stochastic programming model from~\cite{Punnakitikashem2008} by integrating nurse staffing decisions over multiple units into the NPA problem.
Moreover, several recent papers consider patient appointment planning in outpatient chemotherapy clinics while simultaneously assigning nurses to patients or taking constraints on nurse availability into account~\cite{Liang2016,Hesaraki2019,Hesaraki2020,Bouras2021,Heshmat2018}.
Schmidt et al.~\cite{Schmidt2013} integrate patient appointment planning into a patient-to-ward assignment problem, where they assume that rooms in the same ward are equal.
They propose a binary integer program for this problem and compare exact and heuristic solution approaches.
Ceschia and Schaerf~\cite{Ceschia2016} consider an integrated planning of PRA with operating room constraints where they allow the postponement of patient appointments, for which they also provide an instance generator and a set of benchmark instances.

%The nurse staffing problem (assigning an appropriate number of nurses to each unit within a ward given a nurse roster) and the nurse-to-patient assignment problem in a neonatal intensive care unit are considered in an integrated fashion by Pesant~\cite{Pesant2016}, who solves constraint programming models for the two problems consecutively in order to obtain good integrated solutions.

% \smallskip

% One general, high-level hyper-heuristic approach that can be used for both the patient-to-bed assignment problem (as defined by Demeester et al. (2020)) and the nurse rostering problem:
% Bilgin et al. (2012)~\cite{Bilgin2012}

% \smallskip

% Patient-bed-assignment with constraints on nurses:
% Thomas et al. (2013)~\cite{Thomas2013} (patient-to-bed assignment with constraints on the required nurse-to-patient ratio)

% \smallskip

% Patient appointments (in outpatient chemotherapy) planned with consideration of nurses:
% Hesaraki et al. (2019) and (2020)~\cite{Hesaraki2019,Hesaraki2020}, Bouras et al. (2021)~\cite{Bouras2021}, Heshmat et al. (2018)~\cite{Heshmat2018} %(https://doi.org/10.1016/j.cie.2018.07.033 - no access)

% Papers that need not be mentioned:
%Lin et al. (2020)~\cite{Lin2020}, Naderi et al. (2021)~\cite{Naderi2021} (do not deal with hospitals)

% \input{OurContribution}
\section{Problem definition} \label{sec:probDef}

In this section, we formally introduce the \emph{integrated patient-to-room and nurse-to-patient assignment (IPRNPA)} problem as well as the sets and parameters that are used to represent the input of the problem.

\medskip

The IPRNPA problem integrates PRA and NPA on the ward level. Thus, the problem consists of assigning patients to rooms and nurses to patients on a hospital ward over a given planning period (typically one or several weeks). This section describes the static version of the problem where, similar to the static version of the PRA problem~\cite{Demeester2010}, all information about the patients (in particular, each patient's admission and discharge times) are known at the beginning of the planning period. The dynamic version of the problem differs from the static version in that information about new patients only becomes known either when they are admitted or a fixed time span before admission, which is analogous to the existing literature on dynamic PRA~\cite{Ceschia2011,Ceschia2012}.

\medskip

We are given a set~$\setP$ of patients, a set~$\setN$ of nurses, and a set~$\setR$ of rooms. Moreover, there exists a (usually small) set~$\setA$ of additional rooms (such as the nursing station) with $\setR\cap\setA=\emptyset$. These additional rooms cannot be used for assigning patients and will only be relevant when computing the walking distances of the nurses. % The set of patients is partitioned into the subsets~$\setF$ of female patients and~$\setM$ of male patients, i.e., $\setP=\setF \dot{\cup} \setM$.

\medskip

% The considered planning period consists of a set~$\setD=\{1,\dots,D\}$ of days, where each day consists of~$k$ shifts (usually $k=3$). The set of all shifts is denoted by $\setS=\{1,\dots,S=3k\}$ and the subset of all night shifts by~$\setNight:=\{s\in \setS\mid s \mod k = 0\}\subseteq\setS$. Both days and shifts are numbered chronologically, so the first $k$~shifts belong to day~$1$ etc.

The considered planning period consists of a set~$\setS=\{1,\dots,S\}$ of shifts, which is partitioned into the subsets~$\setEarly$ of early shifts, $\setLate$ of late shifts, and $\setNight$ of night shifts. The shifts are numbered chronologically starting with an early shift and ending with a night shift. Hence, the first early, late, and night shift are numbered $1$, $2$, and, $3$, respectively, and belong to the first day, whereas the second early, late, and night shift belong to the second day and so on.

% \medskip

% Each room~$r\in\setR$ has a shift-independent number of available beds denoted by~$\numbeds(r)$ (usually between~$1$ and~$4$) that defines the maximum number of patients that can be assigned to room~$r$ during any single shift. Moreover, rooms can feature different types of equipment (which is assumed to be uniform across all beds in the room) that might be required for certain patients. The set of possible equipment types is denoted by~$\setE$ and the parameter~$\eqpres(r,e)\in\{0,1\}$ specifies whether equipment type~$e$ is present in room~$r$ ($\eqpres(r,e)=1$) or not ($\eqpres(r,e)=0$).

\medskip

A feasible PRA demands that each patient~$p\in\setP$ is assigned to exactly one room~$r\in\setR$ during each shift between their admission shift $\adshift(p)\in\setEarly$ and their discharge shift $\dishift(p)\in\setNight$, which denotes the first and the last shift, respectively, of the patient's stay on the ward. In particular, this means that patients are always admitted and discharged in the morning between a night shift and the following early shift, as is the case in many real hospital wards. If patient~$p$ has already been on the ward during the last shift of the previous planning period, the value~$\adshift(p)$ is set to~$0\notin\setS$, and if patient~$p$ will still be on the ward after the last shift~$S$ of the current planning period, the value~$\dishift(p)$ is set to~$S+1\notin\setS$. Patient transfers between rooms are possible and are assumed to take place at most once a day for each patient between a night shift and an early shift. Transfers are, however, undesirable for both patients and nurses, and should, thus, be minimized. If patient~$p$ has already been on the ward during the last shift of the previous planning period, the room~$y^{\prev}(p)\in\setR$ that the patient has been assigned to during this shift is also given as an input. This allows the evaluation of transfers that happen between the last (night) shift of the previous planning period and the first (early) shift of the current planning period.
%Moreover, there is a patient-specific upper bound~$\maxtransfers(p)$ on the number of times a single patient~$p\in\setP$ can be transferred during their stay (which might depend, for instance, on the number of shifts the patient stays on the ward). If a patient has already been on the ward before the start of the current planning period, we assume that transfers in previous planning periods have already been subtracted from this upper bound beforehand, so the provided value corresponds to the remaining allowed transfers.

\smallskip

Requirements concerning the PRA include respecting the capacity of each room $r\in\setR$, which is given shift-independent number of available beds denoted by~$\numbeds(r)$ (usually between~$1$ and~$4$) that defines the maximum number of patients that can be assigned to room~$r$ during any single shift. Moreover, depending on their specific condition, a patient might benefit from certain types of equipment in their room during certain shifts, so, during these shifts, they should be assigned to a room that features this type of equipment if possible. The set of possible equipment types is denoted by~$\setE$. The types of equipment that are present in room~$r$ are represented by the subset~$\setE(r)\subseteq\setE$, and the shift-specific types of desired equipment of patient~$p$ during shift~$s$ are represented by the subset~$\setE(p,s)\subseteq\setE$. %The shift-specific parameter~$\eqreq(p,e,s)\in\{0,1\}$ indicates whether patient~$p$ requires equipment type~$e$ during shift~$s$ ($\eqreq(p,e,s)=1$) or not ($\eqreq(p,e,s)=1$) and the parameter~$\eqpres(r,e)\in\{0,1\}$ specifies whether equipment type~$e$ is present in room~$r$ ($\eqpres(r,e)=1$) or not ($\eqpres(r,e)=0$).
Additionally, gender-mixed rooms should be avoided if possible. To this end, the set of patients is partitioned into the subsets~$\setF$ of female patients and~$\setM$ of male patients (i.e., $\setP=\setF \dot{\cup} \setM$) and the number of gender-mixed rooms should be minimized across all shifts. %it is required that no room can be assigned both female and male patients during any shift.
Finally, each patient~$p\in\setP$ has an associated age group computed as $\agegroup(p) = \lfloor\frac{\age(p)}{10}\rfloor$ with $\age(p)$ denoting the age of the patient in years. Age groups are relevant since large age differences between patients who are simultaneously assigned to the same room are known to result in inconvenience for the patients and should, thus, be avoided.

% \medskip

% Each patient~$p\in\setP$ has an age (in years) denoted by $\age(p)$, which is relevant since large age differences between patients that are simultaneously assigned to the same room are known to result in inconvenience for the patients and should, thus, be avoided. Moreover, each patient~$p\in\setP$ has an associated admission shift $\adshift(p)\in\setS$ and discharge shift $\dishift(p)\in\setS$, which denote the first and the last shift, respectively, in which the patient needs to be assigned to a room and a nurse. If patient~$p$ has already been on the ward during the last shift of the previous planning period, the value~$\adshift(p)$ is set to~$0\notin\setS$, and if patient~$p$ will still be on the ward after the last shift~$S$ of the current planning period, the value~$\dishift(p)$ is set to~$S+1\notin\setS$.

\medskip

Concerning the NPA part of the problem, the nurse roster for the planning period is given as an input. Here, for each nurse~$n\in\setN$, we are given the subset~$\setS(n)\subseteq\setS$ of shifts that the nurse is assigned to. A feasible NPA should assign each patient~$p\in\setP$ to exactly one nurse~$n\in\setN$ with $s\in\setS(n)$ during each shift~$s\in\setS$ between the patient's admission shift and discharge shift. Since nurses work at most one shift per day, any patient staying on the ward for at least two shifts must necessarily be assigned to different nurses during different shifts. In order to improve continuity of care, however, the number of different nurses who treat a single patient should be minimized. %and a patient-specific upper bound~$\maxdiffnurses(p)$ on the number of nurses patient~$p\in\setP$ is assigned to should be respected.

\smallskip

Further requirements on the NPA include respecting nurse skill level requirements of the patients. Each nurse~$n\in\setN$ has a skill level~$\skilllevel(n)$ and each patient~$p\in\setP$ requires a certain minimum skill level~$\skillreq(p,s)$ during each shift~$s\in\setS\setminus\setNight$ between their admission shift and their discharge shift. The set of possible skill levels is denoted by $\setL=\{1,2,3\}$, where $1=\textnormal{trainee}$, $2=\textnormal{regular}$, and $3=\textnormal{experienced}$. While an experienced nurse (skill level~$3$) can take care of any patient, the assignment should ensure that regular nurses (skill level~$2$) and trainees (skill level~$1$) are only assigned patients whose required skill level during a specific shift does not exceed the nurse's skill level.

\smallskip

Moreover, patients can induce different, shift-dependent workloads for nurses, and these workloads are to be distributed fairly among the nurses - both during each single shift and overall. The workload resulting from taking care of patient~$p\in\setP$ during shift~$s\in\setS$ is expressed by a nonnegative number~$\wload(p,s)$ and depends on the age group of the patient, on their specific condition, on the time since admission, and on whether the shift is a day shift or a night shift. A fair distribution of workload among the nurses is then achieved by defining a maximum desired workload~$\maxload(n,s)$ for each nurse~$n\in\setN$ during each shift~$s\in\setS(n)$ and ensuring that this maximum workload is not exceeded if possible while at the same time making sure that the assigned workloads relative to the respective maxima do not differ too much between nurses during single shifts as well as overall.

\medskip

Important considerations that are influenced by both the PRA as well as the NPA involve trying to assign all patients in a room to the same nurse during each shift (which minimizes the spread of infections across rooms by nurses and also leads to patients being disturbed less by other nurses entering their room) as well as the minimization of the walking distances of nurses between different rooms. For the evaluation of the walking distances, we are given a nonnegative number~$\dist(r,r')$ for each two rooms~$r,r'\in\setR$ that specifies the walking distance between rooms~$r$ and~$r'$. Similarly, for each additional room~$a\in\setA$ and each regular room~$r\in\setR$, there is a nonnegative number~$\dist(a,r)$ that specifies the walking distance between additional room~$a$ and room~$r$. The actual walking distances of nurses during a shift~$s\in\setS$ then depend on the patients they are assigned to during shift~$s$ and on how frequently nurses are expected to walk in a \emph{circular pattern} directly from patient to patient during shift~$s$ (e.g., when rounds are made during early shifts) and on how frequently they are expected to walk in a \emph{star-like pattern} directly from additional rooms (such as the nursing station) to patients and back (e.g., when a patient calls for a nurse). These expected (absolute) frequencies are represented by two nonnegative parameters~$\wpat^{\circ}(s)$ (\emph{circular}) and $\wpat^{\star}(s)$ (\emph{star-like}), respectively.

\section{Mixed integer programming formulation}\label{sec:MIP}

In order to model the IPRNPA problem introduced in Section~\ref{sec:probDef} mathematically, we now present a formulation of the problem as a MIP. %As a point of comparison, we additionally consider a natural decomposition of this MIP into the two interacting subproblems, which is described in Section~\ref{subsec:decomposition}.

\medskip

\noindent
The following lists summarize the sets and parameters introduced in the previous section as well as the decision variables used in the MIP:
%The MIP formulation uses the following sets, parameters, and decision variables:

\medskip

\noindent
\textbf{Sets:}
\smallskip
\begin{itemize}[align=left,labelindent=0pt,labelwidth=2.7cm,labelsep*=1em,leftmargin=!]
\item[$\setP$\quad] set of patients (index~$p$)
\item[$\setF$\quad] subset of female patients
\item[$\setM$\quad] subset of male patients
\item[$\setN$\quad] set of nurses (index~$n$)
\item[$\setNprev(p)$] subset of nurses that patient~$p\in\setP$ has already been assigned to during at least one shift in a previous planning period
\item[$\setR$\quad] set of rooms (index~$r$)
\item[$\setA$\quad] set of additional rooms (e.g., the nursing station) (index~$a$)
\item[$\setS=\{1,\dots,S\}$\quad] set of shifts (index~$s$)
\item[$\setEarly$\quad] subset of early shifts
\item[$\setLate$\quad] subset of late shifts
\item[$\setNight$\quad] subset of night shifts
\item[$\setS(n)$] subset of shifts that nurse~$n\in\setN$ is assigned to
\item[$\setE$\quad] set of possible equipment types in the rooms (index~$e$)
\item[$\setE(r)$\quad] subset of equipment types present in room~$r\in\setR$ %(index~$e$)
\item[$\setE(p,s)$\quad] subset of desired equipment types of patient~$p\in\setP$ during shift~$s\in\setS$
% \item[$\setD=\{1,\dots,D\}$\quad] set of days of the planning period (index~$d$)
\item[$\setL=\{1,2,3\}$\quad] set of possible skill levels of nurses (index~$l$), where $1=\textnormal{trainee}$, $2=\textnormal{regular}$, $3=\textnormal{experienced}$
\end{itemize}

\medskip

\noindent
\textbf{Parameters:}
\smallskip
\begin{itemize}[align=left,labelindent=0pt,labelwidth=2.7cm,labelsep*=1em,leftmargin=!]
\item[$\adshift(p)$\quad] shift~$s\in\setEarly$ during which patient~$p\in\setP$ is admitted (first shift in which a bed is required for patient~$p$). The value is set to~$0$ if patient~$p$ has already been on the ward during the last shift of the previous planning period
\item[$\dishift(p)$\quad] shift~$s\in\setLate$ during which patient~$p\in\setP$ is discharged (last shift in which a bed is required for patient~$p$). The value is set to~$S+1$ if patient~$p$ will still be on the ward after the last shift~$S$ of the planning period
\item[$y^{\prev}(p)$\quad] room~$r\in\setR$ that patient~$p\in\setP$ with $\adshift(p)=0$ has been assigned to during the last shift of the previous planning period %The value is set to ``none'' if patient~$p$ has not been on the ward during the last shift of the previous planning period.
\item[$\numbeds(r)$\quad] nonnegative integer (most likely from $\{1,2,3,4\}$) specifying the number of beds in room~$r\in\setR$
% \color{gray}
% \item[$\maxtransfers(p)$] nonnegative integer specifying the maximum allowed number of transfers of patient~$p\in\setP$
% \color{black}
%\item[$\prevtransfers(p)$\quad] number of transfers of patient~$p\in\setP$ in previous planning periods
\item[$\agegroup(p)$] age group of patient~$p\in\setP$ computed as $\agegroup(p) = \lfloor\frac{\age(p)}{10}\rfloor$ with $\age(p)$ denoting the age of the patient in years
\item[$\skilllevel(n)$\quad] skill level of nurse~$n\in\setN$ (possible values are $1,2,3$)
\item[$\skillreq(p,s)$\quad] minimum skill level of a nurse required by patient~$p\in\setP$ during shift~$s\in\setS\setminus\setNight$ (possible values are $0,1,2$)
% \item[$\maxdiffnurses(p)$] nonnegative integer specifying the maximum number of different nurses that patient~$p\in\setP$ can be assigned to during their stay
%\item[$\assign(n,s)$\quad] $1$ if nurse~$n\in\setN$ is assigned to shift~$s\in\setS$,\\ $0$ otherwise
% \item[$\num(s)$\quad] number of nurses assigned to shift~$s\in\setS$
\item[$\wload(p,s)$\quad] nonnegative number specifying the workload resulting from taking care of patient~$p\in\setP$ during shift~$s\in\setS$
\item[$\maxload(n,s)$\quad] nonnegative number specifying the maximum workload allowed for nurse~$n\in\setN$ during shift~$s\in\setS(n)$
% \item[$\fairnessfactor$\quad] number $\geq1$ specifying the maximum factor by which the average relative workloads (relative to the allowed maxima) of any two nurses~$n,n'\in\setN$ are allowed to differ 
% \item[$\fairnessfactor(s)$\quad] number $\geq1$ specifying the maximum factor by which the relative workloads (relative to the allowed maxima) of any two nurses~$n,n'\in\setN$ assigned to a shift~$s\in\setS$ are allowed to differ during this shift
%\item[$\gender(p)$\quad] $w$~if patient~$p\in\setP$ is female,\\ $m$ if patient~$p\in\setP$ is male
%\item[$\eqpres(r,e)$\quad] $1$ if room~$r\in\setR$ has equipment~$e\in\setE$,\\ $0$ otherwise
%\item[$\eqreq(p,e,s)$\quad] $1$ if patient~$p\in\setP$ requires a room having equipment~$e\in\setE$ during shift~$s\in\setS$,\\ $0$ otherwise
\item[$\dist(r,r')$\quad] nonnegative number specifying the walking distance between rooms~$r,r'\in\setR$ (where $\dist(r,r')=\dist(r',r)$ for all~$r,r'\in\setR$, i.e., distances are symmetric) %resulting for a nurse being assigned patients in both rooms during shift~$s\in\setS$ (where $\dist(r,r',s)=\dist(r',r,s)$ for all~$r,r'\in\setR$, $s\in\setS$, i.e., distances are symmetric)
\item[$\dist(a,r)$\quad] nonnegative number specifying the walking distance between additional room~$a\in\setA$ and room~$r\in\setR$
% \item[$\dist(a,r,s)$\quad] nonnegative number specifying the walking distance resulting for a nurse between additional room~$a\in\setA$ and room~$r\in\setR$ when the nurse is assigned a patient in room~$r$ during shift~$s\in\setS$ (where $\dist(a,r,s)=\dist(r,a,s)$ for all~$a\in\setA$, $r\in\setR$, $s\in\setS$)
% \item[$\fit(p,n,s)$\quad] number (negative, zero, or positive) specifying the quality of fit of a patient~$p\in\setP$ and a nurse~$n\in\setN$ during shift~$s\in\setS(n)$ (where zero means average fit). In particular, includes the nurse's preference for taking care of the patient.
\item[$\wpat^{\circ}(s)$\quad] nonnegative weight for different walking patterns depending on the shift~$s\in\setS$. A high value of $\wpat^\circ(s)$ indicates that most nurses walk directly from patient to patient during shift~$s$ (circular pattern)
\item[$\wpat^{\star}(s)$\quad] nonnegative weight for different walking patterns depending on the shift~$s\in\setS$. A high value of $\wpat^\star(s)$ indicates that most nurses walk directly from additional rooms such as the nursing station to patient rooms and back during shift~$s$ (star-like pattern)
% \item[$\wpat^{\circ}(s)$, $\wpat^{\star}(s)$\quad] nonnegative weights for different walking patterns depending on the shift~$s\in\setS$. A high value of $\wpat^\circ(s)$ indicates that most nurses walk directly from patient to patient during shift~$s$, whereas a high value of $\wpat^\star(s)$ indicates that most nurses walk directly from additional rooms such as the break room to patient rooms and back during shift~$s$.
%\item[$\prevassigned(p,n)$\quad]$1$ if patient~$p\in\setP$ has already been assigned to nurse~$n\in\setN$ during at least one shift in a previous planning period,\\ $0$ otherwise
\end{itemize}

\medskip

\noindent
\textbf{Decision variables:}
\smallskip
\begin{itemize}[align=left,labelindent=0pt,labelwidth=2.7cm,labelsep*=1em,leftmargin=!]
\item[$y_{p,r,s}$\quad] binary variable indicating whether patient~$p\in\setP$ is assigned to room~$r\in\setR$ during early shift~$s\in\setEarly$ \\ (only defined if $\adshift(p)\leq s\leq \dishift(p)$, i.e., if patient~$p$ is on the ward during early shift~$s$)
\item[$\femaleinroom_{r,s}$\quad] binary variable indicating whether at least one female patient is assigned to room~$r\in\setR$ during early shift $s\in\setEarly$
\item[$\maleinroom_{r,s}$\quad] binary variable indicating whether at least one male patient is assigned to room~$r\in\setR$ during early shift~$s\in\setEarly$
\item[$\vio^{\viogender}_{r,s}$\quad] binary variable indicating whether more than one gender is accommodated in room~$r\in\setR$ during shift~$s\in\setEarly$
\item[$\trans_{p,s}$\quad] binary variable indicating whether patient~$p\in\setP$ is transferred to a different room after night shift~$s\in(\setNight\cup\{0\})\setminus \{S\}$ (and before early shift~$s+1$) \\
(only defined if $\adshift(p)\leq s\leq\dishift(p)-1$)
% \item[$y_{p,r,s}$\quad] $1$ if patient~$p\in\setP$ is assigned to room~$r\in\setR$ during early shift~$s\in\setEarly$, $0$ otherwise \\ (only defined if $\adshift(p)\leq s\leq \dishift(p)$, i.e., if patient~$p$ is on the ward during early shift~$s$)
% \item[$\femaleinroom_{r,s}$\quad] $1$ if at least one female patient is assigned to room~$r\in\setR$ during shift~$s\in\setEarly$, $0$ otherwise
% \item[$\maleinroom_{r,s}$\quad] $1$ if at least one male patient is assigned to room~$r\in\setR$ during shift~$s\in\setEarly$, $0$ otherwise \\
% \item[$\vio^{\viogender}_{r,s}$\quad] $1$ if in room~$r\in\setR$ more than one gender is accommodated during shift~$s\in\setEarly$, $0$ otherwise\\
% \item[$\trans_{p,s}$\quad] $1$ if patient~$p\in\setP$ is transferred to a different room after night shift~$s\in(\setNight\cup\{0\})\setminus \{S\}$ (and before early shift~$s+1$), $0$ otherwise
% (only defined if $\adshift(p)\leq s\leq\dishift(p)-1$)
\item[$\agegroup^{\mathrm{max}}_{r,s}$] nonnegative fractional variable representing the maximum age group among all patients~$p\in\setP$ assigned to room~$r\in\setR$ during early shift~$s\in\setEarly$
\item[$\agegroup^{\mathrm{min}}_{r,s}$] nonnegative fractional variable representing the minimum age group among all patients~$p\in\setP$ assigned to room~$r\in\setR$ during early shift~$s\in\setEarly$
%
% Variables of second model
%
\item[$x_{p,n,s}$\quad] binary variable indicating whether patient~$p\in\setP$ is assigned to nurse~$n\in\setN$ during shift~$s\in\setS$,\\ (only defined if $s\in\setS(n)$ and $\adshift(p)\leq s\leq \dishift(p)$, i.e., if if nurse~$n$ is assigned to shift~$s$ and patient~$p$ is on the ward during shift~$s$)
\item[$\vio^{\skill}_{p,s}$\quad] binary variable indicating whether patient~$p\in\setP$ is assigned to a nurse with a lower skill level than required during shift~$s\in\setS\setminus\setNight$\\ (only defined if $\adshift(p)\leq s\leq\dishift(p)$ and $\skillreq(p,s)\geq2$, i.e., if patient~$p$ is on the ward during shift~$s$ and requires at least an experienced nurse during shift~$s$)
%%% The \vio^{\skill}_{p,s} are actually implemented as NONNEGATIVE FRACTIONAL variables in the model since their values will be 0 or 1 anyway by constraint\ref{constr:nurse-skills}
\item[$\everassigned_{p,n}$\;] binary variable indicating whether patient~$p\in\setP$ is assigned to nurse~$n\in\setN$ during at least one shift~$s\in\setS$
% \item[$x_{p,n,s}$\quad] $1$ if patient~$p\in\setP$ is assigned to nurse~$n\in\setN$ during shift~$s\in\setS$,\\ $0$ otherwise\\ (only defined if $s\in\setS(n)$ and $\adshift(p)\leq s\leq \dishift(p)$, i.e., if if nurse~$n$ is assigned to shift~$s$ and patient~$p$ is on the ward during shift~$s$)
% \item[$\vio^{\skill}_{p,s}$\quad] $1$ if patient~$p\in\setP$ is assigned to a nurse with a lower skill level than required during shift~$s\in\setS\setminus\setNight$,\\ $0$ otherwise\\ (only defined if $\adshift(p)\leq s\leq\dishift(p)$ and $\skillreq(p,s)\geq2$, i.e., if patient~$p$ is on the ward during shift~$s$ and requires at least an experienced nurse during shift~$s$)
% %%% The \vio^{\skill}_{p,s} are actually implemented as NONNEGATIVE FRACTIONAL variables in the model since their values will be 0 or 1 anyway by constraint\ref{constr:nurse-skills}
% \item[$\everassigned_{p,n}$\quad] $1$ if patient~$p\in\setP$ is assigned to nurse~$n\in\setN$ during at least one shift~$s\in\setS$,\\ $0$ otherwise
\item[$\vio^{\load}_{n,s}$\quad] nonnegative fractional variable representing the excess load assigned to nurse~$n\in\setN$ during shift~$s\in\setS$\\ (only defined if $s\in\setS(n)$, i.e., if nurse~$n$ is assigned to shift~$s$)
% \item[$\vio^{\fair}_{n,n',s}$\quad] $1$ if the relative workloads of nurses~$n,n'\in\setN$ during shift~$s\in\setS$ differ by more than a factor of $\fairnessfactor(s)$,\\ $0$ otherwise\\ (only defined if $s\in\setS(n)\cap\setS(n')$, i.e., if nurses~$n,n'$ are both assigned to shift~$s$)
% \item[$\vio^{\fair}_{n,n'}$\quad] $1$ if the average relative workloads of nurses~$n,n'\in\setN$ differ by more than a factor of $\fairnessfactor$,\\ $0$ otherwise
\item[$\vio^{\fair}_{n,n',s}$\quad] nonnegative fractional variable representing the excess in relative workload (relative to the desired maximum) of nurse~$n\in\setN$ compared to nurse~$n'\in\setN$ during shift~$s\in\setS$
\item[$\vio^{\fair}_{n,n'}$\quad] nonnegative fractional variable representing the overall excess in relative workload (relative to the desired maximum) of nurse~$n\in\setN$ compared to nurse~$n'\in\setN$
\item[$\inroom_{n,r,s}$\quad] binary variable indicating whether nurse~$n\in\setN$ is assigned at least one patient in room~$r\in\setR$ during shift~$s\in\setS$\\ (only defined if $s\in\setS(n)$, i.e., if nurse~$n$ is assigned to shift~$s$)
% \item[$\inroom_{n,r,s}$\quad] $1$ if nurse~$n\in\setN$ is assigned at least one patient in room~$r\in\setR$ during shift~$s\in\setS$,\\ $0$ otherwise\\ (only defined if $s\in\setS(n)$, i.e., if nurse~$n$ is assigned to shift~$s$)
% \color{gray}
% \item[$\onenurse_{r,s}$\quad] $1$ if all patients in room~~$r\in\setR$ are assigned to the same nurse during shift~$s\in\setS$,\\ $0$ otherwise
% \color{black}
\item[$\dist_{n,s}$\quad] nonnegative fractional variable representing the total walking distance for nurse~$n\in\setN$ during shift~$s\in\setS$ (only defined if $s\in\setS(n)$, i.e., if nurse~$n$ is assigned to shift~$s$)
% \item[$\dist^{\circ}_{n,s}$\quad] nonnegative fractional variable representing the circular walking distance for nurse~$n\in\setN$ during shift~$s\in\setS$ (only defined if $s\in\setS(n)$, i.e., if nurse~$n$ is assigned to shift~$s$)
\item[$\bothrooms_{n,r,r',s}$] binary variable indicating whether nurse~$n\in\setN$ is assigned patients in both room~$r\in\setR$ and room~$r'\in\setR$ during shift~$s\in\setS$\\ (only defined if $\assign(n,s)=1$, i.e., if nurse~$n$ is assigned to shift~$s$)
\end{itemize}

\medskip

The objective function of the MIP that is to be minimized consists of a weighted sum of several separate objectives. These objectives include those classically considered in the PRA problem (objectives~\ref{objective:patient_transfers}--\ref{objective:equipmentviolation}) and the NPA problem (objectives~\ref{objective:continuity_of_care}--\ref{objective:soft_constraints}). Moreover, two additional objectives are considered that rely explicitly on the interaction of the two problems: Objective~\ref{objective:one-nurse} considers the assigning the minimum number of nurses per room during each shift, while objective~\ref{objective:walking_distance} minimizes the walking distances of nurses. The weights of the different objectives in the weighted sum are chosen based on the existing literature and discussions with our partner hospital (see Section~\ref{sec:exp_results} for the concrete values used in our computational experiments).\\

% \noindent
% \textbf{Quality of fit objective}
% \begin{enumerate}[resume,label=(\arabic*)]
% \item Maximization of the quality of fit of patients and nurses they are assigned to across all patients, nurses, and shifts (could also be weighted differently for different patients and / or nurses and / or shifts): \label{objective:qualtity_of_fit}
% \begin{align*}
% \max \sum_{\begin{array}{c}\scriptsize p\in\setP,n\in\setN,s\in\setS(n):\\ \scriptsize \adshift(p)\leq s\leq\dishift(p)\end{array}} x_{p,n,s}\cdot\fit(p,n,s)
% \end{align*}
% \end{enumerate}
\noindent
\textbf{Patient transfers objective}
\smallskip
\begin{enumerate}[resume,label=(\arabic*)]
\item Minimization of the number of patient transfers across all patients and shifts (could also be weighted differently for different patients and / or different shifts): \label{objective:patient_transfers}
\begin{align*}
\min \sum_{\substack{p\in\setP,s\in(\setNight\cup\{0\})\setminus\{S\}:\\ \adshift(p)\leq s\leq\dishift(p)-1}}\trans_{p,s}
\end{align*}
\end{enumerate}
\textbf{Patient inconvenience objective}
\smallskip
\begin{enumerate}[resume,label=(\arabic*)]
\item Minimization of the age group difference across all rooms and shifts: \label{objective:inconvenience}
\begin{align*}
\min \sum\limits_{r \in \setR, s \in \setEarly} (\agegroup^{\mathrm{max}}_{r,s} - \agegroup^{\mathrm{min}}_{r,s})
\end{align*}
\end{enumerate}
\clearpage\noindent
\textbf{Gender mixing objective}
\smallskip
\begin{enumerate}[resume,label=(\arabic*)]
\item Minimization of gender mixing across all rooms and shifts: \label{objective:gendermixing}
\begin{align*}
\min \sum\limits_{r \in \setR, s \in \setEarly} \vio^{\viogender}_{r,s} 
\end{align*}
\end{enumerate}
\textbf{Equipment violation objective}
\smallskip
\begin{enumerate}[resume,label=(\arabic*)]
\item Minimization of required equipment violation across all rooms and shifts: \label{objective:equipmentviolation}
\begin{align*}
\min \sum\limits_{\substack{p \in \setP, r \in \setR, s \in \setEarly:\\ \adshift(p)\leq s\leq \dishift(p) \\ \text{and } \setE(p,s)\setminus\setE(r)\neq \emptyset}} y_{p,r,s} 
\end{align*}
\end{enumerate}
\textbf{Continuity of care objective}
\smallskip
\begin{enumerate}[resume,label=(\arabic*)]
\item Minimization of the number of different nurses that treat each patient across all patients (could also be weighted differently for different patients): \label{objective:continuity_of_care}
\begin{align*}
\min \sum_{p\in\setP,n\in\setN\setminus \setNprev(p)} \everassigned_{p,n}
\end{align*}
% \begin{align*}
% \min \sum_{p\in\setP}\left(\sum_{n\in\setN} \everassigned_{p,n}-\prevassigned(p,n)\right)
% \end{align*}
\end{enumerate}
\textbf{Penalization of skill level requirements and undesired workload distributions objective}
\smallskip
\begin{enumerate}[resume,label=(\arabic*)]
\item Minimization of violations of skill level requirements of patients and undesired workload distributions for nurses: %(weights for the different kinds of violation variables could be introduced here to model that not all violations are equally bad):
\label{objective:soft_constraints}
\begin{align*}
\min &\sum_{\substack{ p\in\setP,s\in\setS\setminus\setNight:\\  \adshift(p)\leq s\leq\dishift(p)
\\\text{and } \skillreq(p,s)\geq2}}\vio^{\skill}_{p,s} +\sum_{n\in\setN,s\in\setS(n)}\vio^{\load}_{n,s}\\ &+\sum_{n,n'\in\setN}\vio^{\fair}_{n,n'}
+\sum_{\substack{\scriptsize n,n'\in\setN,\\ s\in\setS(n)\cap\setS(n')}}\vio^{\fair}_{n,n',s}
\end{align*}
\end{enumerate}
% \color{gray}
% \textbf{Assigning one nurse to patients in the same room objective}
% \begin{enumerate}[resume,label=(\arabic*)]
% \item Maximization of the number of room-shift pairs in which all patients in the room are assigned to the same nurse during the shift: \label{objective:one-nurse}
% \begin{align*}
% \max &\sum\limits_{r \in R, s \in S} \onenurse_{r,s}
% \end{align*}
% \end{enumerate}
\textbf{Assigning the minimum number of nurses per room objective}
\smallskip
\begin{enumerate}[resume,label=(\arabic*)]
\item Minimization of the number of nurses assigned to rooms across all shifts: \label{objective:one-nurse}
\begin{align*}
\min &\sum\limits_{n \in N, r \in R, s \in S} \inroom_{n,r,s}
\end{align*}
\end{enumerate}
\textbf{Walking distances objective}
\smallskip
\begin{enumerate}[resume,label=(\arabic*)]
\item Minimization of the walking distances across all nurses and shifts %(could also be weighted differently for different nurses and / or different shifts):
\label{objective:walking_distance}
\begin{align*}
\min \sum_{n\in\setN,s\in\setS(n)} \dist_{n,s}
\end{align*}
\end{enumerate}

\bigskip

\noindent
The constraints of the MIP can be formulated as follows:\\
\medskip

\noindent
(I)~\textbf{Assignment of patients to rooms}
\smallskip
\begin{enumerate}[resume,label=(\arabic*)]
\item Each patient~$p\in\setP$ is assigned to exactly one room~$r\in\setR$ during each early shift~$s\in\setEarly$ between their admission and discharge: %(and to no room at all in all the remaining shifts):
\label{constr:patient-room-assignment}
\begin{align*}
\sum_{r\in\setR} y_{p,r,s} = 1 & \quad\forall p\in\setP,\, s\in\setEarly: \adshift(p)\leq s\leq\dishift(p) %\\
%\sum_{r\in\setR} y_{p,r,s} = 0 & \quad\forall p\in\setP,\, s\in\setS: s\leq\adshift(p)-1 \text{ or } s\geq \dishift(p)+1
\end{align*}
% \item Patients are only assigned to rooms that have all the required equipment: \label{constr:room-features}
% \begin{align*}
% y_{p,r,s} = 0 \quad\forall p\in\setP,\, r\in\setR,\, s\in\setEarly: & \adshift(p)\leq s\leq \dishift(p),\, \setE(p,s)\setminus\setE(r)\neq \emptyset
% \end{align*}
\item No room~$r\in\setR$ can be assigned more than $\numbeds(r)$ patients during any early shift~$s\in\setEarly$: \label{constr:room-capacities}
\begin{align*}
\sum_{p\in\setP: \adshift(p)\leq s\leq \dishift(p)} y_{p,r,s} \leq \numbeds(r) & \quad\forall r\in\setR,\, s\in\setEarly
\end{align*}
\item The variable~$\femaleinroom_{r,s}$ ($\maleinroom_{r,s}$) is set to one if at least one female (male) patient is assigned to room~$r$ during early shift~$s\in\setEarly$: \label{constr:set_in_room_vars}
\begin{align*}
    y_{p,r,s} & \leq\femaleinroom_{r,s} & \forall p\in\setF,\, r\in\setR,\, s\in\setEarly: \adshift(p)\leq s\leq \dishift(p)\\
    y_{p,r,s} & \leq\maleinroom_{r,s} & \forall p\in\setM,\, r\in\setR,\, s\in\setEarly: \adshift(p)\leq s\leq \dishift(p)
\end{align*}
\item No room~$r\in\setR$ should be assigned both female and male patients during any early shift~$s\in\setEarly$: \label{constr:no-gender-mix}
\begin{align*}
\femaleinroom_{r,s} + \maleinroom_{r,s} \leq 1 + \vio^{\viogender}_{r,s} & \quad\forall r\in\setR,\, s\in\setEarly
\end{align*}
% \item (Cut) For any two rooms~$r,r'\in\setR$ that have the same size and equipment, the one whose total/average distance to the additional rooms~$a\in\setA$ is smallest (room~$r$) should be used first: \label{cut:room-assignment}
% \begin{align*}
% -\left(\sum_{a\in\setA}\dist(a,r)\right)\cdot(\femaleinroom_{r,s} + \maleinroom_{r,s}) \leq -\left(\sum_{a\in\setA}\dist(a,r')\right)\cdot(\femaleinroom_{r',s} + \maleinroom_{r',s})\\ - M\cdot\sum_{p\in\setP} (2-y_{p,r',s}-y_{p,r',s-1}) \\ \quad\forall s\in\setS\setminus\{1\},\, r,r'\in\setR: \numbeds(r) = \numbeds(r') \text{ and } \setE(r)=\setE(r')
% \end{align*}
\end{enumerate}

\noindent
(II)~\textbf{Patients transfers}
\smallskip
\begin{enumerate}[resume,label=(\arabic*)]
\item The patient transfer variables~$\trans_{p,s}$ are set correctly for each patient~$p\in\setP$ and each night shift~$s\in\setNight$ between their admission and discharge: \label{constr:transfer-vars}
\begin{samepage}
\begin{align*}
y_{p,r,s+1} - y_{p,r,s-2} \leq \trans_{p,s} \quad &\forall p\in\setP,\, r\in\setR, \, s\in\setNight\setminus\{S\}:\\
&\adshift(p)\leq s\leq\dishift(p)-1
\end{align*}
\end{samepage}
\item The patient transfer variables~$\trans_{p,0}$ that indicate a transfer between the last shift of the previous planning period (shift~$0$) and the first (early) shift of the current planning period (shift~$1$) are set correctly for each patient~$p\in\setP$ with $\adshift(p)=0$: \label{constr:transfer-vars2} % Remark: We assume here that the first shift is an early shift
\begin{align*}
y_{p,r,1} \leq \trans_{p,0} & \quad\forall p\in\setP,\, r\in\setR\setminus \{y^{\prev}(p)\}:\adshift(p)=0
\end{align*}
% \color{gray}
% \item Every patient $p\in\setP$ can be transferred at most $\maxtransfers(p)$ times during their stay:
% \begin{align*}
%     \sum_{s\in(\setNight\setminus\{S\})\cup\{0\}:\adshift(p)\leq s\leq\dishift(p)-1} \trans_{p,s} \leq \maxtransfers(p)-\prevtransfers(p) & \quad\forall p\in\setP
% \end{align*}
% \color{red}
\end{enumerate}

\clearpage\noindent
(III)~\textbf{Inconvenience of patients}
\smallskip
\begin{enumerate}[resume,label=(\arabic*)]
%\item The variables $\age^{\mathrm{max}}_{r,s}$ are restricted by the maximum age of prior occupants in a room~$r\in\setR$ on shift~$s\in\setS$: \label{constr:ageMax-1}
%\begin{align*}
%\age^{\mathrm{max}}(r,s) \geq \age^{\mathrm{max}}_{r,s}  &  \quad\forall r\in\setR,\, s\in\setS
%\end{align*}

\item The variable $\agegroup^{\mathrm{max}}_{r,s}$ is restricted by the maximum age group of patients in room~$r\in\setR$ during early shift~$s\in\setEarly$: \label{constr:ageMax}
\begin{align*}
\agegroup^{\mathrm{max}}_{r,s} \geq \agegroup(p) \cdot y_{p,r,s} \quad &\forall p\in\setP,\, r\in\setR,\, s\in\setEarly:\\
& \adshift(p)\leq s \leq \dishift(p)
\end{align*}

%\item The variables $\age^{\mathrm{min}}_{r,s}$ are restricted by the minimum age of prior occupants in a room~$r\in\setR$ on shift~$s\in\setS$: \label{constr:ageMin-1}
%\begin{align*}
%\age^{\mathrm{min}}_{r,s} \leq \age^{\mathrm{min}}(r,s)  &  \quad\forall r\in\setR,\, %s\in\setS
%\end{align*}

\item The variable $\agegroup^{\mathrm{min}}_{r,s}$ is restricted by the minimum age group of patients in room~$r\in\setR$ during early shift~$s\in\setEarly$: \label{constr:ageMin}
\begin{align*}
\agegroup^{\mathrm{min}}_{r,s} \leq \agegroup(p) + 12 \cdot (1 - y_{p,r,s}) \quad &\forall p \in \setP,\, r\in\setR,\, s\in\setEarly:\\
& \adshift(p)\leq s \leq \dishift(p)
\end{align*}
Note that the coefficient~$12$ on the right-hand side needs to be increased if patients with ages~$130$ or older (age group~$13$ or higher) are present.
\smallskip

\item The variable $\agegroup^{\mathrm{min}}_{r,s}$ is set to zero if no patients are assigned to room~$r\in\setR$ during early shift~$s\in\setEarly$: \label{constr:ageMin-2}
\begin{align*}
\agegroup^{\mathrm{min}}_{r,s} \leq 12\cdot \sum_{p\in\setP: \adshift(p)\leq s \leq \dishift(p)} y_{p,r,s} &  \quad\forall r\in\setR,\, s\in\setEarly
\end{align*}
Note that the coefficient~$12$ on the right-hand side needs to be increased if patients with ages~$130$ or older (age group~$13$ or higher) are present.
\smallskip

\item The value of the variable $\agegroup^{\mathrm{max}}_{r,s}$ must not be smaller than the value of the variable $\agegroup^{\mathrm{min}}_{r,s}$ for any room~$r\in\setR$ and any early shift~$s\in\setEarly$:
\label{constr:ageMinMax}
\begin{align*}
\agegroup^{\mathrm{min}}_{r,s} \leq \agegroup^{\mathrm{max}}_{r,s} &  \quad\forall r\in\setR,\, s\in\setEarly
\end{align*}
Note that this constraint is not required for the correctness of the model and is only used to improve solution times.
\end{enumerate}

% Constraints of second model
% \color{black}
\bigskip
\noindent
(IV)~\textbf{Assignment of patients to nurses}
\smallskip
\begin{enumerate}[resume,label=(\arabic*)]
\item Each patient~$p\in\setP$ is assigned to exactly one nurse~$n\in\setN$ with $s\in\setS(n)$ during each shift~$s\in\setS$ between their admission and discharge: %(and to no nurse in all the remaining shifts):
\label{constr:patient-nurse-assignment}
\begin{align*}
\sum_{n\in\setN:s\in\setS(n)} x_{p,n,s} = 1 & \quad\forall p\in\setP,\, s\in\setS: \adshift(p)\leq s\leq\dishift(p) %\\
%\sum_{n\in\setN} x_{p,n,s} = 0 & \quad\forall p\in\setP,\, s\in\setS: s\leq\adshift(p)-1 \text{ or } s\geq \dishift(p)+1
\end{align*}
%\item During a shift, patients are only assigned to nurses that have been assigned to this shift in the nurse roster: \label{constr:nurse-roster}
%\begin{align*}
%x_{p,n,s} \leq \assign(n,s) & \quad\forall p\in\setP,\, n\in\setN,\, s\in\setS
%\end{align*}
\clearpage\noindent
\item For each patient~$p\in\setP$ and each shift~$s\in\setS\setminus\setNight$ with $\skillreq(p,s)\geq 2$, the variable~$\vio^{\skill}_{p,s}$ is set to one if the patient is not assigned to a nurse with the required skill level: \label{constr:nurse-skills}
\begin{align*}
\sum_{\substack{n\in\setN:s\in\setS(n) \text{ and}\\ \skilllevel(n)\geq \skillreq(p,s)}} x_{p,n,s} = 1 - \vio^{\skill}_{p,s} \quad &\forall p\in\setP,\, s\in\setS\setminus\setNight:\\ & \adshift(p)\leq s\leq \dishift(p)\\ & \text{and } \skillreq(p,s)\geq2
\end{align*}
% \begin{align*}
% \sum_{\scriptsize\begin{array}{c} n\in\setN:s\in\setS(n)\\ \scriptsize \text{and } \skilllevel(n)\geq \skillreq(p,s)\end{array}} & x_{p,n,s} = 1 - \vio^{\skill}_{p,s} \\ \quad \forall p\in\setP,\, s\in\setS\setminus\setNight:\, & \adshift(p)\leq s\leq \dishift(p) \text{ and } \skillreq(p,s)\geq2
% \end{align*}
% Old version of the constraint (where this has been a hard constraint)
% \item Patients are only assigned to nurses that have all the required skills: \label{constr:nurse-skills}
% \begin{align*}
% x_{p,n,s} = 0 & \quad\forall p\in\setP,\, n\in\setN,\, s\in\setS: \skillreq(p,k,s)=1 \text{ and } \skillpres(n,k)=0 \text{ for some skill}~k
% \end{align*}
\item The variable~$\everassigned_{p,n}$ is set to one if and only if $n\in\setNprev(p)$ or patient~$p\in\setP$ is assigned to nurse~$n\in\setN$ during at least one shift: \label{constr:ever_assigned}
\begin{align*}
x_{p,n,s} \leq \everassigned_{p,n} & \quad\forall p\in\setP,\, n\in\setN,\, s\in\setS: s\in\setS(n)\\
& \quad\text{and} \adshift(p)\leq s \leq \dishift(p)\\
\everassigned_{p,n} = 1 & \quad\forall p\in\setP,\, n\in\setNprev(p)
\end{align*}
\begin{align*}
\everassigned_{p,n} \leq \sum_{\substack{s\in\setS: s\in\setS(n) \text{ and}\\ \adshift(p)\leq s \leq \dishift(p)}} \hspace{-4mm}x_{p,n,s} \quad\forall p\in\setP,\, n\in\setN\setminus\setNprev(p)
\end{align*}
% \item Each patient~$p$ is allowed to be assigned to at most $\maxdiffnurses(p)$ different nurses during their stay: \label{constr:ub_diff_nurses}
% \begin{align*}
% \sum_{n\in\setN}\everassigned_{p,n} \leq \maxdiffnurses(p) & \quad\forall p\in\setP
% \end{align*}
\end{enumerate}

\noindent
(V)~\textbf{Workload of nurses}
\smallskip
\begin{enumerate}[resume,label=(\arabic*)]
\item For each nurse~$n\in\setN$ and each shift~$s\in\setS(n)$, any workload exceeding $\maxload(n,s)$ leads to a corresponding increase of the variable~$\vio^{\load}_{n,s}$: \label{constr:nurse-load}
\begin{align*}
\sum_{p\in\setP: \adshift(p)\leq s\leq \dishift(p)} x_{p,n,s}\cdot\wload(p,s) \leq\; & \maxload(n,s) + \vio^{\load}_{n,s} \\  &\forall n\in\setN,\, s\in\setS(n)
\end{align*}
% \item Fair distribution of workload during each shift: If the relative workloads (relative to the allowed maxima) of any two nurses~$n,n'\in\setN$ assigned to a shift~$s\in\setS$ differ by more than a factor of $\fairnessfactor(s)$ during this shift, the variable~$\vio^{\fair}_{n,n',s}$ is set to one: \label{constr:fairness_nurse-load_shift}
% \begin{align*}
% & \; \sum_{p\in\setP: \adshift(p)\leq s\leq \dishift(p)} x_{p,n,s}\cdot\frac{\wload(p,s)}{\maxload(n,s)} \\
% \leq & \; \fairnessfactor(s)\cdot\sum_{p\in\setP: \adshift(p)\leq s\leq \dishift(p)} x_{p,n',s}\cdot\frac{\wload(p,s)}{\maxload(n',s)}+M\cdot\vio^{\fair}_{n,n',s}\\ & \forall n,n'\in\setN,\, s\in\setS(n)\cap\setS(n')
% \end{align*}
% Here, the factor~$M$ (which represents a sufficiently large integer) is only needed since the actual loads of a nurse could be larger than the allowed maxima if the soft constraints~\ref{constr:nurse-load} are violated.
\item Fair distribution of workload during each shift: For any two nurses $n,n'\in\setN$ assigned to a shift~$s\in\setS$, if the relative workload (relative to the desired maximum) of nurse~$n$ during shift~$s$ exceeds the relative workload of nurse~$n'$ during shift~$s$, the
variable~$\vio^{\fair}_{n,n',s}$ must be increased correspondingly:\label{constr:fairness_nurse-load_shift}
\begin{align*}
&\sum_{p\in\setP: \adshift(p)\leq s\leq \dishift(p)} x_{p,n,s}\cdot\frac{\wload(p,s)}{\maxload(n,s)} \\
\leq& \sum_{p\in\setP: \adshift(p)\leq s\leq \dishift(p)} x_{p,n',s}\cdot\frac{\wload(p,s)}{\maxload(n',s)}+\vio^{\fair}_{n,n',s} \\
& \forall n,n'\in\setN,\, s\in\setS(n)\cap\setS(n')
\end{align*}
% \item Fair distribution of workload overall: If the average relative workloads (relative to the allowed maxima) of any two nurses~$n,n'\in\setN$ differ by more than a factor of $\fairnessfactor$, the variable~$\vio^{\fair}_{n,n'}$ is set to one: \label{constr:fairness_nurse-load}
% \begin{align*}
% & \; \sum_{s\in\setS(n)}\sum_{p\in\setP: \adshift(p)\leq s\leq \dishift(p)} x_{p,n,s}\cdot\frac{\wload(p,s)}{\maxload(n,s)} \\
% \leq & \; \fairnessfactor\cdot\sum_{s\in\setS(n')}\sum_{p\in\setP: \adshift(p)\leq s\leq \dishift(p)} x_{p,n',s}\cdot\frac{\wload(p,s)}{\maxload(n',s)}+M\cdot\vio^{\fair}_{n,n'}\quad\forall n,n'\in\setN
% \end{align*}
% Here, the factor~$M$ (which represents a sufficiently large integer) is only needed since the actual loads of a nurse could be larger than the allowed maxima if the soft constraints~\ref{constr:nurse-load} are violated.
\item Fair distribution of workload overall: For any two nurses~$n,n'\in\setN$, if the relative workload (relative to the desired maximum) of nurse~$n$ exceeds the relative workload of nurse~$n'$, the
variable~$\vio^{\fair}_{n,n'}$ must be increased correspondingly:\label{constr:fairness_nurse-load}
\begin{align*}
&\sum_{s\in\setS(n)}\sum_{p\in\setP: \adshift(p)\leq s\leq \dishift(p)} x_{p,n,s}\cdot\frac{\wload(p,s)}{\maxload(n,s)} \\
\leq& \sum_{s\in\setS(n')}\sum_{p\in\setP: \adshift(p)\leq s\leq \dishift(p)} x_{p,n',s}\cdot\frac{\wload(p,s)}{\maxload(n',s)}+\vio^{\fair}_{n,n'}\\ & \forall n,n'\in\setN
\end{align*}
\end{enumerate}

\noindent
(VI)~\textbf{Assignment of all patients in the same room to the same nurse}
\smallskip
\begin{enumerate}[resume,label=(\arabic*)]
\item The variables~$\inroom_{n,r,s}$ are set correctly for each nurse~$n\in\setN$, each room~$r\in\setR$, and each shift~$s\in\setS(n)$: \label{constr:inroom}
\begin{align*}
\inroom_{n,r,s} & \geq x_{p,n,s} + y_{p,r,s} - 1 &\forall p\in\setP,\, n\in\setN,\, r\in\setR,\, s\in\setS(n)\cap\setEarly:\\ && \adshift(p)\leq s \leq \dishift(p) \\
\inroom_{n,r,s} & \geq x_{p,n,s} + y_{p,r,s-1} - 1 &\forall p\in\setP,\, n\in\setN,\, r\in\setR,\, s\in\setS(n)\cap\setLate:\\ && \adshift(p)\leq s \leq \dishift(p) \\
\inroom_{n,r,s} & \geq x_{p,n,s} + y_{p,r,s-2} - 1 &\forall p\in\setP,\, n\in\setN,\, r\in\setR,\, s\in\setS(n)\cap\setNight:\\ && \adshift(p)\leq s \leq \dishift(p)
\end{align*}
% \color{gray}
% \item The variable~$\onenurse_{r,s}$ is set to zero if any two patients in room~$r\in\setR$ are assigned to different nurses during shift~$s\in\setS$: \label{constr:oneNurse-1}
% \begin{align*}
%  \inroom_{n,r,s} + \inroom_{n',r,s} - 1 \leq N \cdot (1 - \onenurse_{r,s}) & \quad\forall r\in\setR,\, n,n'\in\setN,\, s\in\setS(n)\cap\setS(n'): n > n'
% \end{align*}
% \item The variable~$\onenurse_{r,s}$ is set to zero if no patients are assigned to room~$r\in\setR$ during shift~$s\in\setS$: \label{constr:oneNurse-2}
% \begin{align*}
% \onenurse_{r,s} & \leq \sum_{p\in\setP:\adshift(p)\leq s \leq \dishift(p)} y_{p,r,s} &  \quad\forall r\in\setR,\, s\in\setEarly \\
% \onenurse_{r,s} & \leq \sum_{p\in\setP:\adshift(p)\leq s \leq \dishift(p)} y_{p,r,s-1} &  \quad\forall r\in\setR,\, s\in\setLate \\
% \onenurse_{r,s} & \leq \sum_{p\in\setP:\adshift(p)\leq s \leq \dishift(p)} y_{p,r,s-2} &  \quad\forall r\in\setR,\, s\in\setNight
% \end{align*}
% \color{black}
\end{enumerate}

\noindent
(VII)~\textbf{Walking distances of nurses}
\smallskip
\begin{enumerate}[resume,label=(\arabic*)]
\item The variables~$\bothrooms_{n,r,r',s}$ are set correctly for each nurse~$n\in\setN$, each two rooms~$r,r'\in\setR$, and each shift~$s\in\setS(n)$:\label{constr:bothrooms}
\begin{align*}
\bothrooms_{n,r,r',s} \geq \inroom_{n,r,s} + \inroom_{n,r',s} - 1\; \forall n\in\setN,\, r,r'\in\setR,\, s\in\setS(n)
\end{align*}
% \item For each room~$r\in\setR$, the work load due to patients assigned to this room is determined for each nurse~$n\in\setN$ and each shift~$s\in\setS$ with $\assign(n,s)=1$: \label{constr:room-load}
% \begin{align*}
% \roomload_{n,r,s} = \sum_{p\in\setP:\adshift(p)\leq s\leq\dishift(p)}\wload(p,s)\cdot y_{p,r,s}\\ &\forall n\in\setN,\, s\in\setS:\assign(n,s)=1
% \end{align*}
% \item The walking distance variables~$\dist^{\circ}_{n,s}$ are set correctly for each nurse~$n\in\setN$ and each shift~$s\in\setS(n)$: \label{constr:walk-dist1}
% \begin{align*}
% \dist^{\circ}_{n,s} \geq \sum_{r\in\setR}\sum_{r'\in\setR} (\inroom_{n,r,s}+\inroom_{n,r',s}-1)\cdot\dist(r,r',s)\quad &\forall n\in\setN,\, s\in\setS(n)
% \end{align*}
% \begin{align*}
% \color{red}\dist^{\circ}_{n,s} \geq \frac{1}{2}\cdot\sum_{r\in\setR}\sum_{r'\in\setR\setminus\{r\}} (\inroom_{n,r,s}+\inroom_{n,r',s}-1)\cdot\dist(r,r',s)\quad &\color{red}\forall n\in\setN,\, s\in\setS(n)
% \end{align*}
% Incorrect version without variables $\bothrooms_{n,r,r',s}$:
% \item The walking distance variables~$\dist_{n,s}$ are set correctly for each nurse~$n\in\setN$ and each shift~$s\in\setS(n)$: \label{constr:walk-dist2}
%  \begin{align*}
%  \dist_{n,s} &= \wpat^\circ(s)\cdot\dist^{\circ}_{n,s} +\wpat^\star(s)\cdot\sum_{a\in\setA}\sum_{r\in\setR} \inroom_{n,r,s}\cdot\dist(a,r,s) &\forall n\in\setN,\, s\in\setS(n)
%  \end{align*}
\item The walking distance variables~$\dist_{n,s}$ are set correctly for each nurse $n\in\setN$ and each shift~$s\in\setS(n)$:\label{constr:walk-dist}
\begin{align*}
\dist_{n,s} &= \wpat^\circ(s)\cdot\frac{1}{2}\cdot\sum_{r\in\setR}\sum_{r'\in\setR\setminus\{r\}} \bothrooms_{n,r,r',s}\cdot\dist(r,r')\\ &+\wpat^\star(s)\cdot\sum_{a\in\setA}\sum_{r\in\setR} \inroom_{n,r,s}\cdot\dist(a,r)\;\, \forall n\in\setN,\, s\in\setS(n)
\end{align*}
\end{enumerate}
\section{Solution methods}\label{sec:solution-methods}

Addressing PRA and NPA challenges in hospitals reveals complexities that hinder timely optimal solutions. The PRA problem, even in isolation, is known to be \textsf{NP}-hard and only unrealistically small instances can be solved to optimality in reasonable time. Thus, it not surprising that the integrated problem is very difficult to solve on instances of a realistic size and solving the MIP provided in the previous section takes a prohibitive amount of time on such instances (see also Section~\ref{sec:exp_results}). Therefore, we now present two different methods for generating good solutions in reasonable computation times. The first method, which is mainly presented as a point of comparison to the completely integrated MIP from Section~\ref{sec:MIP}, is a sequential approach based on a natural decomposition of the integrated MIP into models for the two interacting subproblems. The second method is an efficient greedy heuristic for the integrated problem, which also allows easy adaptions to various dynamic problem versions.

\subsection{Sequential solution approach}\label{subsec:decomposition}

As a point of comparison to the completely integrated MIP presented in the previous section and to potentially achieve faster computation times, we consider the sequential solution of the two submodels for the PRA part and the NPA part. Here, we first solve the PRA part to determine the PRA and then solve the NPA part to determine the NPA based on the given PRA. Note that this approach of solving the two assignment problems sequentially will of course not yield optimal solutions for the integrated problem in general.

\medskip

The PRA part of the model consists of the patient transfers objective~\ref{objective:patient_transfers}, the patient inconvenience objective~\ref{objective:inconvenience}, the gender mixing objective~\ref{objective:gendermixing}, and the equipment violations objective~\ref{objective:equipmentviolation} as well as the constraints concerning the assignment of patients to rooms~(I), the patients transfers~(II), and the inconvenience of patients~(III) and the associated decision variables.

\medskip

The NPA part of the model contains the continuity of care objective~\ref{objective:continuity_of_care}, the penalization of skill level requirements and undesired workload distributions objective~\ref{objective:soft_constraints}, the objective for assigning the minimum number of nurses per room~\ref{objective:one-nurse}, and the walking distance objective~\ref{objective:walking_distance}. As constraints, the NPA part contains those concerning the assignment of patients to nurses~(IV), the workload of nurses~(V), the assignment of all patients in the same room to the same nurse~(VI), and the walking distance of nurses~(VII). Concerning variables, the NPA model contains all variables appearing in these objectives and constraints, but the~$y_{p,r,s}$ are not decision variables anymore since their values are carried over from the solution obtained for the PRA part of the model.

\subsection{Heuristic solution approach}\label{subsec:heuristic}

Due to the computationally challenging nature of the PRA problem, heuristic solution methods are often used to tackle this problem in the recent literature~\cite{Guido2018,Schaefer2019}. Naturally, the computational challenges become even greater when considering PRAs and NPAs jointly in an integrated optimization problem as we do in this paper. Therefore, we now introduce a heuristic solution approach for the IPRNPA problem, which extends the heuristic for the PRA problem suggested in~\cite{Schaefer2019}. The presented heuristic balances computational effort and solution quality with the goal of obtaining a practically applicable procedure for timely, satisfactory outcomes.

\medskip

Our heuristic is initially tailored to the static version of the IPRNPA problem, but its scope extends beyond the static setting. In fact, the heuristic does not require to have all information about patients available in advance and can easily be extended to a variety of dynamic settings. As an illustrative example, instead of restricting patient transfers between rooms to happen between a night shift and the following early shift as in the static IPRNPA problem, the heuristic could handle patient transfers whenever new information about patients is obtained. The same applies for NPAs, which could alternatively be allowed also during shifts whenever new information becomes available. Thus, even complex dynamic problem settings based on dynamic PRA problems such as the one considered by Ceschia and Schaerf~\cite{Ceschia2011,Ceschia2012} can be handled without difficulties.

%In the context of dynamic PRA, the heuristic demonstrates its adaptability by seamlessly transitioning to an online mode, all the while accommodating diverse configurations. One illustrative instance involves the option to execute transfers throughout each iteration of the heuristic, or alternatively, confining them to predetermined time points. The same applies to the assignment of nurses to patients and the possible re-allocations during shifts. Thus, even complex dynamic problem formulations based on dynamic PRA settings such as the one used by Ceschia and Schaerf~\cite{Ceschia2011,Ceschia2012} can be incorporated without difficulties.

\subsubsection{Algorithm description}\label{subsubsec:general-alg}
The idea of the heuristic is to assign patients to rooms and nurses to patients in a greedy fashion. To this end, the heuristic considers the days of the planning period chronologically, where each day is represented by the corresponding early shift. For each day, the heuristic iteratively fixes both the PRA and the NPA for only this day jointly for a single patient in a way that yields the lowest current contribution to the objective function. After each such assignment, the current objective function contributions of the remaining possible assignments for the day are updated before the assignments for the next patient are fixed. Once the room and nurse assignments have been fixed for all patients that are on the ward during the considered day, the heuristic moves on to the next day.
During this iterative process, we have to take into account that decisions for the two kinds of assignments are made on different time scales. Concerning PRAs, a decision about the room the patient is assigned to is made only once per day before the start of the early shift, whereas, concerning NPAs, three different nurse must be assigned to a patient for each day (for the early, late, and night shift) since each nurse works at most one shift per day.

\medskip

%The heuristic treats each day of the planning period separately in chronological order, where each day is represented by the corresponding early shift~$s\in\setEarly$. 
We now describe our heuristic, whose pseudocode is shown in Algorithm~\ref{alg:heuristic}, more formally. Here, we first describe the algorithm without the heterogeneity check between patient admission and discharge times that is represented by the heterogeneity matrix~$\hetmatrix$ in the pseudocode and is motivated and explained afterwards in Section~\ref{subsubsec:heterogeneity-matrix}.

Upon initialization of the heuristic, all relevant assignment variables~$x_{p,n,s}$ and~$y_{p,r,s}$ are set to~$0$ to indicate that no assignments have been made so far. Afterwards, the days of the planning period are considered in chronological order represented by their early shifts.
When considering a day of the planning period represented by the corresponding early shift~$s\in\setEarly$, the set of relevant patients for which room and nurse assignments are to be made for this day is denoted by $\setP(s) \coloneqq\{p\in\setP:\adshift(p)\leq s \leq \dishift(p)\}$. Note that~$\setP(s)$ also includes patients that have been on the ward already on the previous day if they have not been discharged yet since the room and nurse assignments for such patients have so far only been fixed up to the previous day. At each point during the assignment process for the day, the set of rooms that are still available (i.e., not yet fully occupied) is denoted by $\setR(s)\subseteq\setR$. The set of nurses who are on duty on the corresponding day is partitioned according to their assigned shifts into $\setN^{\text{early}}(s)$, $\setN^{\text{late}}(s)$, and $\setN^{\text{night}}(s)$. Consequently, the Cartesian product $\setNcomb(s) \coloneqq \setN^{\text{early}}(s) \times \setN^{\text{late}}(s) \times \setN^{\text{night}}(s)$ corresponds to all possible ordered triples of nurses that can potentially be assigned to a patient during the three shifts of the day that starts with the early shift~$s\in\setEarly$. The potential contributions to the objective function for each triple~$(p,n^{\text{comb}},r)$ consisting of a patient~$p\in \setP(s)$, a triple~$n^{\text{comb}}\in\setNcomb(s)$ of nurses, and an available room~$r\in\setR(s)$ are stored in a contribution table denoted by $\contribution$.

\medskip

Whenever the heuristic starts considering some day of the planning period starting with early shift $s\in\setEarly$, the contribution table is filled with the current contributions that correspond to the potential assignments of patients~$p\in\setP(s)$ to nurse combinations $n^{\text{comb}}\in\setNcomb(s)$ and rooms~$r\in\setR(s)$. These current contributions are computed based on the changes to the objective function value that would currently be induced by the corresponding assignments. Here, the objective function is the same weighted sum of the objectives~\ref{objective:patient_transfers}--\ref{objective:walking_distance} as in the MIP described in Section~\ref{sec:MIP}.

The heuristic then fixes the assignments of rooms and nurse combinations for the day represented by early shift~$s$ for all patients in~$\setP(s)$ in a greedy fashion. This is done by iteratively identifying the triple~$(p,n^{\text{comb}},r)$ that has the lowest contribution value in $\contribution$ and setting the values of the corresponding assignment variables~$x_{p,n,s}$ and~$y_{p,r,s}$ to~$1$, which means that patient~$p$ is assigned to the nurse combination~$n^{\text{comb}}$ and the room~$r$. 
Subsequently, the allocated patient~$p$ is removed from~$\setP(s)$ and the room~$r$ is removed from~$\setR(s)$ in case that the assignment has resulted in full occupancy of room~$r$ on the corresponding day. Moreover, the contribution table $\contribution$ needs to be updated by removing all entries corresponding to patient~$p$ and possibly room~$r$. In addition, the update involves adjusting the objective contributions for the remaining possible assignments of the day to ensure that they accurately represent the new state of PRA and NPAs. For example, the previously-calculated contribution to the gender mixing objective for room~$r$ is updated for all remaining patients in~$\setP(s)$ with a different gender than the patient~$p$ just assigned. After the updated contribution table has been computed, the heuristic continues with the next iteration based on the updated contribution table until $\setP(s)=\emptyset$ (i.e., all assignments for the current day have been fixed) and the procedure continues with the next day.

\medskip

The heuristic terminates once the last day of the planning period has been considered, i.e., once all PRA and NPAs have been fixed for the whole planning period.

\subsubsection{Heterogeneity check between patients}\label{subsubsec:heterogeneity-matrix}
A point that is not taken into account in the heuristic as described so far is the coordination of arrivals and discharges of patients that are assigned to the same room. If occupancy levels are high, this could lead to the unnecessary creation of gender-mixed rooms or avoidable patient transfers. For example, assigning several male patients with very similar arrival and discharge shifts to different rooms instead of the same room might leave no other rooms available for female patients admitted shortly afterward unless some of the male patients are transferred. In order to avoid such unfavorable incidents explicitly, we next describe a new patient-to-patient heterogeneity measure, whose values are stored in a heterogeneity matrix. This heterogeneity matrix is used in the heuristic to foster the assignment of patients who are %admitted and 
discharged at similar times to the same room by explicitly taking %admission and
discharge shifts into account during the computation of possible objective function contributions.

\medskip

% The proposed patient-to-patient heterogeneity measure is based on the difference between two patients' admission shifts as well as on the difference between their discharge shifts. Formally, for two patients~$p,p'\in\setP$, we let $\text{ad\_shift\_diff}(p,p')$ denote the absolute difference between the admission shifts $\adshift(p)$ and $\adshift(p')$, while $\text{di\_shift\_diff}(p,p')$ denotes the absolute difference between the discharge shifts $\dishift(p)$ and $\dishift(p')$. The heterogeneity for the two patients is then calculated by the following formula\footnote{In case that $\adshift(p)=\adshift(p')$ and $\dishift(p)=\dishift(p')$, where the argument of the logarithm equals zero, we set $\text{het}(p,p')$ to zero. In all other cases, the value is strictly positive.}:
% \begin{align*}
% \text{het}(p,p') \coloneqq \ln{\sqrt{\text{ad\_shift\_diff}(p,p')^2 + \text{di\_shift\_diff}(p,p')^2}}
% \end{align*}
The proposed patient-to-patient heterogeneity measure is based on the difference between two patients' %admission shifts as well as on the difference between their
discharge shifts. Formally, for two patients~$p,p'\in\setP$, we let 
%$\text{ad\_shift\_diff}(p,p')$ denote the absolute difference between the admission shifts $\adshift(p)$ and $\adshift(p')$, while 
$\text{di\_shift\_diff}(p,p')$ denote the absolute difference between the discharge shifts $\dishift(p)$ and $\dishift(p')$. The heterogeneity for the two patients is then calculated by the following formula\footnote{In case that %$\adshift(p)=\adshift(p')$ and 
$\dishift(p)=\dishift(p')$, where the argument of the logarithm equals zero, we set $\text{het}(p,p')$ to zero. In all other cases, the value is strictly positive.}:
\begin{align*}
\text{het}(p,p') \coloneqq \ln{ \text{di\_shift\_diff}(p,p')}
\end{align*}
% This means that the heterogeneity value~$\text{het}(p,p')$ is the natural logarithm of the Euclidean distance of the two pairs $(\adshift(p),\dishift(p))$ and $(\adshift(p'),\dishift(p'))$ when interpreted as points in the plane. Here, the logarithm is used in order to limit the growth of the values in cases where the planning period is long, where very large values would otherwise occur. Note that, since the logarithm is a strictly increasing function, a lower heterogeneity value for two patients signifies more similar admission and discharge dates and, thus, a better fit between the two patients for being assigned to the same room.
This means that the heterogeneity value~$\text{het}(p,p')$ is the natural logarithm of the %Euclidean 
difference between the $\dishift(p)$ and $\dishift(p')$. Here, the logarithm is used in order to limit the growth of the values in cases where the planning period is long, where %very 
large values would otherwise occur. Note that, since the logarithm is a strictly increasing function, a lower heterogeneity value for two patients signifies more similar discharge dates and, thus, a better fit between the two patients for being assigned to the same room.

\medskip

The heterogeneity values of all patient pairs are calculated upon initialization of the heuristic and stored in a heterogeneity matrix denoted by $\hetmatrix$. Here, since the heterogeneity values are symmetric (i.e., $\text{het}(p,p')=\text{het}(p',p)$ for all~$p,p'\in\setP$), it suffices to compute the values above the main diagonal of the matrix in order to improve efficiency. The heterogeneity values of patient pairs are then used whenever the current objective contributions are computed for an assignment that involves assigning a patient to a room to which other patients have already been assigned on the same day. In this case, the appropriately weighted maximum of the heterogeneity values between the new patient and the already assigned patients is added as an additional summand to the current objective contributions.

\alglanguage{pseudocode}
\renewcommand\algorithmicrequire{\textbf{Input:}}
\renewcommand\algorithmicensure{\textbf{Output:}}
\begin{algorithm}
    \caption{IPRNPA heuristic}\label{alg:heuristic}
    \begin{footnotesize}
        \begin{algorithmic}[1]
            \Require An instance of the IPRNPA problem
            \Ensure patient-to-nurse assignments~$x$, patient-to-room assignments~$y$
            \State $x_{p,n,s} \gets 0$ for all $p\in\setP, n\in\setN, s\in \setS$
            \State $y_{p,r,s} \gets 0$ for all $p\in\setP, r\in\setR, s\in \setEarly$
            \State $\hetmatrix \gets \mathrm{calculateHeterogeneityMatrix}$
            \For{$s\in\setEarly$ in increasing order}
                \State $\setP(s) \gets \{p\in\setP:\adshift(p)\leq s \leq \dishift(p)\}$
                \State $\setNcomb(s) \gets \setN^{\text{early}}(s)\times \setN^{\text{late}}(s)\times \setN^{\text{night}}(s)$
                \State $\setR(s)\gets \setR$
                \State $\contribution \gets \mathrm{calculateContributionTable}(\setP(s), \setNcomb, \setR(s),\hetmatrix)$
                \While {$\setP(s) \neq \emptyset$}
                    \State $(p, n^{\text{comb}}, r) \gets \mathrm{argmin}(\contribution)$
                    \State $y_{p,r,s} \gets 1$
                    \State $(n^{\text{early}}, n^{\text{late}}, n^{\text{night}}) \gets n^{\text{comb}}$
                    \State $x_{p,n^{\text{early}},s}, x_{p,n^{\text{late}},s}, x_{p,n^{\text{night}},s} \gets 1$
                    \State $\setP(s) \gets \setP(s)\setminus \{p\}$
                    \If{room~$r$ is fully occupied}
                        \State $\setR(s) \gets \setR(s)\setminus \{r\}$
                    \EndIf
                    \State $\contribution \gets \mathrm{updateContributionTable}(\contribution, p, n^{\text{comb}}, r, \hetmatrix)$
                \EndWhile
                \EndFor
                %\State $\mathrm{return}(x, y)$
        \end{algorithmic}
    \end{footnotesize}
\end{algorithm}

\section{Instance generator and case study}\label{sec:instances}
%\section{Data and test instances}\label{sec:instances}

We now present a detailed description of the parameterized instance generator that we developed for creating realistic test instances of the IPRNPA problem. Afterwards, we present the structure and key data of the real-world instances obtained from our partner hospital. The concrete parameter values used for generating the artificial instances as well as the numbers of instances considered will be described later in Section~\ref{sec:exp_results}.

\subsection{Instance generator} \label{sec:InstanceGenerator}
In order to generate larger numbers of realistic test instances, we developed a parameterized instance generator for the IPRNPA problem. The source code of the instance generator is publicly available on GitHub at \url{https://github.com/TabeaBrandt/instance_generation_integrated_beds_and_staff_planning.git}.

\medskip

The instance generator offers the possibility to create a specifiable number of test instances based on user-defined parameters. These parameters include the number of instances to be created, as well as the option to specify the length of the planning period in weeks. The remaining input parameters can be divided into two main categories: room-related parameters and nurse-related parameters. In the context of room-related parameters, aspects such as the number of patient rooms, their capacity (single, double, triple, or quadruple rooms), occupancy rate, presence of additional rooms, and possible types of room equipment play a central role. In the context of nurse-related parameters, significant options for fine-tuning are available. This includes setting the maximum desired workloads for nurses, as well as specifying the skill levels to be considered. Concerning the total number of nurses that are available and the nurse roster, the generator operates in two modes: manual and automatic. In the manual mode, the user specifies the number of nurses explicitly, based on which the generator strives to create a feasible nurse roster using the binary integer program~(BIP) outlined in Appendix~\ref{sec:A1-nurse-rostering}. In the automatic mode, on the other hand, the generator increases the number of nurses until a feasible nurse roster can be generated using the BIP. %With these parameters the generator generates a realistic and feasible nurse roster by using the binary integer program outlined in Appendix~\ref{sec:A1-nurse-rostering}.

\medskip

The output of the generator is a set of instances of the specified cardinality. In each instance, a defined number of rooms is described, which is divided into single, double, triple, and quadruple rooms according to the specified distribution of room sizes. The rooms are assigned randomly selected equipment from the specified types of possible equipment. Additional rooms such as nursing stations or storage rooms can be added, and a minimum of one such room is required to calculate walking distances based on the star-like walking pattern. The weighting factors that determine the importance of the circular versus the star-like walking pattern depend on the type of shift. The circular pattern is favored during early shifts, a more equal split is favored during late shifts, and the star-like pattern is favored during night shifts.

\medskip

The number of nurses is derived from the number of patient rooms, the distribution of the number of beds per room, and the number of nurses of each skill level that is required per shift, all of which are specified when creating an instance. In instances where three skill levels are specified, we assume that 20\% are experienced nurses (skill level~3), 60\% are regular nurses (skill level~2) and that 20\% are trainees (skill level~1). For instances with two skill levels, we assumed that 80\% are skill level~2 and 20\% are skill level~1. Additionally, the maximum desired workload associated with each skill level is set to 10 for skill level~1, 12.5 for skill level~2, and 15 for skill level~3. This represents an average nurse-to-patient ratio of 1:4, 1:5, and 1:6, respectively, during each shift.
%a ratio of nurses to patients is assumed for each shift 1:8 for the early and late shifts and 1:12 for the night shift. The max working load for each nurse is a shift specified as an instance generator parameter.

% In order to derive a feasible nurse roster for our instances, we built a simple nurse rostering model based on the description presented in the first International Nurse Rostering Competition (INRC) 2010~\cite{Haspeslagh2014}. As the focus of this work is not on the roster, we do not explain it in detail here but provide our used model in the appendix.
% Similar to INRC, we determine the roster for the planning period considering one ward. We use a subset of the described constraints to get a simple, yet still realistic nurse roster. These constraints include that the number of required nurses per shift must be met, not more than a maximum number of shifts can be assigned to a nurse within the time horizon and unwanted shift patterns for nurses are prevented.

\medskip

The patients are generated based on the room configuration and the desired occupancy level. Each patient is assigned a ten-year age group uniformly sampled from 20--30 to 90--100, and, an admission shift, based on the number of rooms and the rooms' capacities. A patient's discharge shift is set as the minimum of the admission shift plus a sampled length of stay (LOS) in days, drawn from a discrete uniform distribution on $\{1,\dots,5\}$, or the last shift of the planning period. Gender is currently assigned based on a 50-50 female-male split. The nurse skill level required for each shift of a patient's LOS is assumed to decrease monotonously. The workload generated by a patient~$p$ during each shift is based on a gamma distribution with $\alpha = 3$, $\beta = 0.5 + \agegroup(p) / 10$, with a minimum of~$1$, a maximum of~$5$, and an exponential smoothing parameter of~$0.1$ that describes the monotonous decrease. The equipment required by the patient is sampled from the types of possible room equipment and assumes monotonously decreasing requirements over the shifts of a patient's LOS. %Recall that equipment is allocated to rooms and is not portable.

\subsection{Real-world instances}\label{sec:real-instances}
To investigate the potential of our methods using real-world data, we also test them on real-world instances from a Short Stay Unit of our partner hospital. As with our instance generator, these instances are publicly available on GitHub at \url{https://github.com/TabeaBrandt/instance_generation_integrated_beds_and_staff_planning.git}.

\medskip

The considered ward is not restricted to a single medical specialty, but only patients requiring care that can be delivered according to a strict protocol are admitted. Consequently, there are no acute admissions to this ward and admission and discharge dates of patients as well as their care requirements are known in advance. The ward is closed on weekends, so each patient's LOS is at most five days.

\medskip

During the week, the nursing staff operates in three shifts (early, late and night) to ensure continuous day and night care. Due to the protocol-based care, the nurses who work on the ward do not specialize in a single medical specialty, but require a broad skill set. There are two nurse skill levels: experienced and trainee. One experienced nurse can take care of four to six patients simultaneously during a shift depending on the patient's care requirements and the nurse's experience, while trainee nurses can take care of about two patients in parallel. There is no particular nurse-to-patient ratio during night shifts, but at least two nurses must always be present. %& The nursing staff operates in three shifts: early, day, and night, ensuring continuous coverage throughout the day and night.

\medskip

The ward consists of 17~patient rooms of varying sizes: four single rooms, 10~double rooms, two triple rooms, and one quadruple room. Therefore, the ward has a total capacity of 34~beds. 
%In reality, however, there are often not enough nurses to staff all beds, so the ward leaves several beds empty.
Additionally, there is one nursing station where the nurses are usually located when they are not attending to patients. We were provided with the floor plan to estimate the walking distances between the rooms, which we calculated according to the shortest walking path between the centers of each room pair.
%Input floor map, or stylized version of it?

\medskip

In our numerical experiments, we use real admission data and nurse rosters of 40~weeks (about nine months), from before the COVID-19 pandemic. %During COVID the ward and/or its nurses were used for COVID care.
Because the ward closes on the weekend, each week in the data can be considered as a separate instance.
We acknowledge that, despite the comprehensive dataset shared by the hospital, certain input data have been omitted due to privacy concerns in order to safeguard individual patient identities.
For example, data on the skill levels required for taking care of individual patients and the resulting patient-specific workloads for nurses have not been provided. These missing data have been generated using the corresponding functions of our instance generator presented in Section~\ref{sec:InstanceGenerator} based on realistic parameter values that have been established in cooperation with Amsterdam University Medical Centers (see Section~\ref{sec:real-instances}).
The data provided led to infeasible instances for two weeks. In one case, a patient was assigned to a shift for which no nurse was on duty, while in the second instance, a shortage in bed capacity during a particular shift rendered it infeasible to accommodate the required patient load.
\section{Experimental results}\label{sec:exp_results}

This section presents our experimental results obtained by testing the MIP as well as the solution methods presented in Section~\ref{sec:solution-methods} on both artificial instances generated by our instance generator from Section~\ref{sec:InstanceGenerator} and the real-world instances described in Section~\ref{sec:real-instances}.

\medskip

All computational experiments were performed on a Linux system running Ubuntu~23.04. The hardware includes an AMD EPYC~7542 processor with 32~CPU cores and 64~threads, operating at a base clock speed of 2.9GHz. The system is equipped with 512~GB of memory. The experiments were implemented using Python~3.11 and Gurobi~10.0.1. To solve the completely integrated MIP presented in Section~\ref{sec:MIP} and the two submodels in the sequential solution approach from Section~\ref{subsec:decomposition}, we dedicated 8~threads to each instance when solving with Gurobi.

\medskip

As described in Section~\ref{sec:MIP}, the objective function to be minimized consists of a weighted sum of several separate objectives.
The specific weights for these objectives were determined based on the existing literature and discussion with our partner hospital. The weights from the existing literature were taken from Demeester et al.~\cite{Demeester2010} and are as follows: the patient transfers objective~\ref{objective:patient_transfers} was assigned a weight of~11, the gender mixing objective~\ref{objective:gendermixing} a weight of~5, and the equipment violation objective~\ref{objective:equipmentviolation} a weight of~5. Furthermore, in consultation with our partner hospital the following weights were established: the patient inconvenience objective~\ref{objective:inconvenience} and the continuity of care objective~\ref{objective:continuity_of_care} were each given a weight of~1, the penalization of skill level requirements and undesired workload distributions objective~\ref{objective:soft_constraints} a weight of~5, the assignment of the minimum number of nurses per room objective~\ref{objective:one-nurse} a weight of~2, and the walking distances objective~\ref{objective:walking_distance} a weight of~0.05. When applying our heuristic, we also considered the heterogeneity values of patient pairs as discussed in Section~\ref{subsubsec:heterogeneity-matrix}. Here, we assigned a weight of~1 to appropriately incorporate this factor.

\subsection{Artificial instances}\label{subsec:results-artificial}
We established a structured framework involving two distinct scenarios, each subdivided into three specific variations and two different lengths of the planning period. This results in a total of 12 scenario-variation-planning period combinations. For each of these 12~combinations, we generated 10~artificial instances using our instance generator described in Section~\ref{sec:InstanceGenerator}.

\medskip

The two scenarios encompass configurations of 30 and 60~beds. Within each scenario, Variation~1 comprises exclusively double rooms, Variation~2 exclusively triple rooms, and Variation~3 encompasses a diverse mix of room types, including single, double, triple, and quadruple rooms. For the 30~beds scenario, this allocation translates to 3~single rooms, 5~double rooms, 3~triple rooms, and 2~quadruple rooms. These numbers are doubled for the 60~beds scenario. Moreover, our investigation includes two different lengths of the planning period for each scenario-variation combination, spanning either~2 or 4~weeks.

\medskip

Throughout our comprehensive analysis, we used input parameters for our instance generator that encompass two distinct equipment types, three possible nurse skill levels, the inclusion of a nursing station as an additional room, and a constant occupancy rate of~85\%.

% \medskip

% As described in Section~\ref{sec:MIP}, the objective function of the MIP that is to be minimized consists of a weighted sum of several separate objectives.
% The specific weights for these objectives were determined based on the existing literature and discussion with our partner hospital (Amsterdam University Medical Centers). The weights from the existing literature were taken from Demeester et al.~\cite{Demeester2010} and are as follows: The patient transfers objective~\ref{objective:patient_transfers} was assigned a weight of~11, the gender mixing objective~\ref{objective:gendermixing} a weight of~5, and the equipment violation objective~\ref{objective:equipmentviolation} a weight of~5. Furthermore, in consultation with our partner hospital the following weights were established: The patient inconvenience objective~\ref{objective:inconvenience} and the continuity of care objective~\ref{objective:continuity_of_care} were each given a weight of~1, the penalization of skill level requirements and undesired workload distributions objective~\ref{objective:soft_constraints} a weight of~5, the assigning the minimum number of nurses per room objective~\ref{objective:one-nurse} a weight of~2, and the walking distances objective~\ref{objective:walking_distance} a weight of~0.05. When applying our heuristic, we also considered the heterogeneity values of patient pairs as discussed in Section~\ref{subsubsec:heterogeneity-matrix}. Here, we assigned a weight of~15 to appropriately incorporate this factor.

\medskip

In addressing the given scenarios, we applied three distinct solution methods to solve the artificial instances, which we will refer to as Methods~1--3 in the following: the MIP from Section~\ref{sec:MIP} (Method~1), the sequential solution approach from Section~\ref{subsec:decomposition} (Method~2) and the heuristic solution approach from Section~\ref{subsec:heuristic} (Method~3). For Method~1, the running time was limited to a maximum of 3~hours (10800 seconds). In the case of Method~2, the same total running time limit was evenly distributed between the two addressed submodels. Additionally, we configured all MIPs to be terminated when the MIP gap falls below~5\%. The results, including objective values and running time data, are presented in Table~\ref{tab:artificalObjectives}.

\begin{table}[ht]
    \centering
    \caption{Artificial instances: Running times and objective values} \label{tab:artificalObjectives}
    \begin{adjustbox}{max width=\linewidth} % Adjust the table width
    \begin{tabular}{cccrrrr}
        \toprule
        \multirow{2}{*}{Scenario} & \multirow{2}{*}{Variation} & \multirow{2}{*}{\shortstack{Planning\\horizon}} & \multicolumn{2}{c}{Running time (seconds)} & \multicolumn{2}{c}{Objective value\tnote{1}} \\
        & & & \multicolumn{1}{c}{Avg} & \multicolumn{1}{c}{Stdev} & \multicolumn{1}{c}{Avg} & \multicolumn{1}{c}{Stdev} \\
        \midrule
        \multicolumn{3}{l}{\textbf{Method 1: MIP}} & & & & \\
\multirow{6}{*}{30 beds} & \multirow{2}{*}{Var. 1 (double rooms)} & 2 weeks & 10800 & 0 & 100\% & 0\% \\
 &  & 4 weeks & 10800 & 0 & 100\% & 0\% \\

 & \multirow{2}{*}{Var. 2 (triple rooms)} & 2 weeks & 10800 & 0 & 100\% & 0\% \\
 &  & 4 weeks & 10800 & 0 & 100\% & 0\% \\

 & \multirow{2}{*}{Var. 3 (mixed rooms)} & 2 weeks & 10800 & 0 & 100\% & 0\% \\
 &  & 4 weeks & 10800 & 0 & 100\% & 0\% \\

\midrule
\multirow{6}{*}{60 beds} & \multirow{2}{*}{Var. 1 (double rooms)} & 2 weeks & 10800 & 0 & 100\% & 0\% \\
 &  & 4 weeks & 10800 & 0 & 100\% & 0\% \\

 & \multirow{2}{*}{Var. 2 (triple rooms)} & 2 weeks & 10800 & 0 & 100\% & 0\% \\
 &  & 4 weeks & 10800 & 0 & 100\% & 0\% \\

 & \multirow{2}{*}{Var. 3 (mixed rooms)} & 2 weeks & 10800 & 0 & 100\% & 0\% \\
 &  & 4 weeks & 10800 & 0 & 100\% & 0\% \\

\midrule
\multicolumn{3}{l}{\textbf{Method 2: Sequential solution approach}} & & & & \\
\multirow{6}{*}{30 beds} & \multirow{2}{*}{Var. 1 (double rooms)} & 2 weeks & 10800 & 0 & 61\% & 2\% \\
 &  & 4 weeks & 10800 & 0 & 46\% & 14\% \\

 & \multirow{2}{*}{Var. 2 (triple rooms)} & 2 weeks & 7829 & 970 & 73\% & 4\% \\
 &  & 4 weeks & 10800 & 0 & 36\% & 10\% \\

 & \multirow{2}{*}{Var. 3 (mixed rooms)} & 2 weeks & 10800 & 0 & 64\% & 2\% \\
 &  & 4 weeks & 10800 & 0 & 38\% & 13\% \\

\midrule
\multirow{6}{*}{60 beds} & \multirow{2}{*}{Var. 1 (double rooms)} & 2 weeks & 10800 & 0 & 25\% & 1\% \\
 &  & 4 weeks & 10800 & 0 & 25\% & 3\% \\

 & \multirow{2}{*}{Var. 2 (triple rooms)} & 2 weeks & 10800 & 0 & 30\% & 23\% \\
 &  & 4 weeks & 10800 & 0 & 24\% & 2\% \\

 & \multirow{2}{*}{Var. 3 (mixed rooms)} & 2 weeks & 10800 & 0 & 24\% & 1\% \\
 &  & 4 weeks & 10800 & 0 & 24\% & 1\% \\

\midrule
\multicolumn{3}{l}{\textbf{Method 3: Heuristic solution approach}} & & & & \\
\multirow{6}{*}{30 beds} & \multirow{2}{*}{Var. 1 (double rooms)} & 2 weeks & 27 & 1 & 88\% & 4\% \\
 &  & 4 weeks & 57 & 1 & 65\% & 19\% \\

 & \multirow{2}{*}{Var. 2 (triple rooms)} & 2 weeks & 38 & 1 & 116\% & 6\% \\
 &  & 4 weeks & 80 & 1 & 57\% & 16\% \\

 & \multirow{2}{*}{Var. 3 (mixed rooms)} & 2 weeks & 30 & 1 & 98\% & 5\% \\
 &  & 4 weeks & 65 & 2 & 58\% & 21\% \\

\midrule
\multirow{6}{*}{60 beds} & \multirow{2}{*}{Var. 1 (double rooms)} & 2 weeks & 888 & 46 & 17\% & 1\% \\
 &  & 4 weeks & 1531 & 62 & 16\% & 0\% \\

 & \multirow{2}{*}{Var. 2 (triple rooms)} & 2 weeks & 1071 & 46 & 27\% & 20\% \\
 &  & 4 weeks & 1784 & 66 & 17\% & 1\% \\

 & \multirow{2}{*}{Var. 3 (mixed rooms)} & 2 weeks & 871 & 45 & 18\% & 1\% \\
 &  & 4 weeks & 1515 & 101 & 17\% & 0\% \\

        \bottomrule
    \end{tabular}
    \end{adjustbox}
    \begin{tablenotes}
        \item \tnote{1} All objective values are provided as proportions (expressed as percentages) of those obtained using Method~1. Values below 100\% indicate superior performance to Method~1.
    \end{tablenotes}
\end{table}

\medskip

The reported results present the objective values as percentages, with the Method~1 objective value serving as the baseline for each instance. These percentages represent the resulting objective proportions, averaged across the 10 instances for each scenario-variation-planning period combination. It is important to emphasize that values below 100\% indicate superior performance to Method~1 since the model aims to minimize the objective function. Notably, Method~2 consistently outperforms Method~1 in all scenario-variation-planning period combinations, while Method~3 solution approach surpasses Method~1 in nearly all instances. When comparing Method~2 to Method~3, it becomes clear that Method~2 excels in the 30 beds scenario but lags behind Method~3 for the 60 beds scenario, except for one case. The standard deviation of the objective values for Method 2 and Method 3 is generally low across most scenario-variation-planning period combinations. However, a few instances with higher standard deviations can be attributed to a handful of strongly deviating instances.

\medskip

When comparing the running times of the three methods, it is evident that Method~1 consistently exceeded the running time limit in all instances, while Method~2 did so in nearly all instances. In contrast, Method~3 displayed remarkable efficiency, requiring, on average, only 27 to 80 seconds for the 30~beds scenario instances and 871 to 1784 seconds for the 60~beds scenarios. This disparity in running times is primarily due to the non-linear growth in nurse combinations for one-day NPAs, as the heuristic assigns three nurses simultaneously. In total, the 30~beds scenario had 21 nurses available to be on duty, while the 60~beds scenario had 31 nurses. Additionally, for each scenario-variation-planning period combination, the running times remained stable with low standard deviation, highlighting the consistency of Method~3's performance. Doubling the planning horizon from 2 to 4 weeks in both scenarios and variations resulted in an approximate doubling of the running time, indicating a linear relationship between running time and planning horizon length.

\medskip

By factoring in the achieved MIP gaps (Table \ref{tab:mipgap_artificial} in appendix) alongside the objective values and running times, we gain insights into the difficulty levels of each scenario-variation-planning period combination. Variation~2, which exclusively considers triple rooms in the 30 beds scenario, emerges as the easiest to solve for both Method~1 and Method~2, while Variations~1 and 3, focusing on double rooms and a mixed set of room types, respectively, seem equally challenging for both methods. Notably, Method~3 exhibits greater efficiency in Variations~1 and 3 compared to Variation 2 when considering running times. MIP gaps exceeding 100\% for Method~1 indicate cases where the initial root relaxation of the MIP model could not be solved within the running time limit, resulting in Gurobi finding suboptimal heuristic solutions. In summary, Method~3 consistently outperforms Method~1 and Method~2 across all cases, excelling in both running time and objective value, particularly in complex scenarios where achieving optimality is notably demanding.

\medskip

Since the considered objective function is a weighted sum of the separate objective functions~\ref{objective:patient_transfers}--\ref{objective:walking_distance}, we also consider the values of these separate objective functions for a more comprehensive analysis. The detailed outcomes are presented in Table~\ref{tab:artificialSearateObjectives}.

\medskip

\begin{table}[ht]
    \centering
    \caption{Artificial instances: Separate objective values differentiated by method} \label{tab:artificialSearateObjectives}
    \begin{adjustbox}{max width=\linewidth} % Adjust the table width
    \begin{tabular}{rlrrrrrr}
        \toprule
        \multicolumn{2}{c}{\multirow{3}{*}{Objective}}  & \multicolumn{2}{c}{\textbf{Method~1}} & \multicolumn{2}{c}{\textbf{Method~2}} & \multicolumn{2}{c}{\textbf{Method~3}} \\
        & & \multicolumn{2}{c}{Obj. value\tnote{1}} & \multicolumn{2}{c}{Obj. value\tnote{1}} & \multicolumn{2}{c}{Obj. value\tnote{1}} \\
        & & \multicolumn{1}{c}{30 beds} & \multicolumn{1}{c}{60 beds} & \multicolumn{1}{c}{30 beds} & \multicolumn{1}{c}{60 beds} & \multicolumn{1}{c}{30 beds} & \multicolumn{1}{c}{60 beds} \\
        \midrule
        \ref{objective:patient_transfers} & Transfers &	100\% & 	100\% & 	0\% & 	21\% & 	1\% & 	3\% \\ 
        \ref{objective:inconvenience} & Inconvenience &	100\% & 	100\% & 	29\% & 	41\% & 	64\% & 	56\% \\ 
        \ref{objective:gendermixing} & Gender mixing & 	100\% & 	100\% & 	2\% & 	2\% & 	40\% & 	28\% \\ 
        \ref{objective:equipmentviolation} & Equipment violation & 	100\% & 	100\% & 	144\% & 	24\% & 	291\% & 	74\% \\ 
        \ref{objective:continuity_of_care} & Continuity of care &	100\% & 	100\% & 	103\% & 	109\% & 	98\% & 	100\% \\ 
        \ref{objective:soft_constraints} & Skill \& workload &	100\% & 	100\% & 	26\% & 	17\% & 	68\% & 	9\% \\ 
        \ref{objective:one-nurse} & Nurses per room &	100\% & 	100\% & 	83\% & 	109\% & 	79\% & 	101\% \\ 
        \ref{objective:walking_distance} & Walking distances &	100\% & 	100\% & 	59\% & 	36\% & 	53\% & 	26\% \\ 
        \bottomrule
    \end{tabular}
    \end{adjustbox}
    \begin{tablenotes}
        \item \tnote{1} All objective values are provided as proportions (expressed as percentages) of those obtained using Method~1 (MIP).
    \end{tablenotes}
\end{table}

Upon examining each objective separately, several noteworthy observations emerge. Method~2 and Method~3 in both scenarios exhibit a notable advantage over Method~1 in strictly minimizing patient transfers (objective~\ref{objective:patient_transfers}). The objectives related to continuity of care (objective~\ref{objective:continuity_of_care}) and the number of nurses per room (objective~\ref{objective:one-nurse}) exhibit similar behavior across Methods~1 to 3. Gender mixing (objective~\ref{objective:gendermixing}) is strongly avoided in Method~2 and to a moderate extent in Method~3. Interestingly, it is worth noting that skill and workload violations (objective~\ref{objective:equipmentviolation}) are more pronounced in Methods~2 and 3 only in the 30~beds scenario, while this issue does not appear as prominently in the 60~beds scenario. These findings underscore each method's distinct strengths and weaknesses in addressing specific objectives.

\subsection{Real-world instances}\label{subsec:results-real-world}
To test all developed methods using real-world data, we sourced data from a Short Stay Unit of our partner hospital as described in Section~\ref{sec:real-instances}. The dataset comprises 40~individual instances, each spanning a planning period of one week.

\medskip

To address the absence of certain input data needed for the IPRNPA problem, we generated the missing information using the corresponding functions implemented in our instance generator introduced in Section~\ref{sec:InstanceGenerator} based on realistic parameter values established in cooperation with our partner hospital. Specifically, we utilized the functions to generate data on nurse skill level requirements of patients and nurse workloads induced by patients. This approach allowed us to create complete instances for our analyses while at the same time maintaining patient privacy and confidentiality.

\medskip

When comparing average (minimum, maximum) parameter values across the considered instances, we observe notable variations. There are approximately~17 (min.~13, max.~23) nurses attending to around~62 (min.~45, max.~76) patients. Patient LOS are around 4.1 (min.~3.6, max.~4.8) shifts, resulting in an average occupancy rate of 61\% (min.~46\%, max.~72\%), while each patient causes a workload of about 2.7 (min.~2.2, max.~3.1) during a day shift, i.e., early and late shift.\footnote{The presented LOS and workload numbers are averages over all patients of an instance.} 
%patient workload and patient length of stay are themselves average values for each instance - here the valuee represent the average/minimum/maximum of the average instance values 
We selected the 20~instances with the highest occupancy rates from the available dataset for further analysis, taking care to exclude holiday times and instances exhibiting unnatural utilization levels. This careful selection process ensures that the chosen instances accurately represent the most demanding and meaningful scenarios for in-depth examination.

\medskip

In tackling the real-world instances, we again applied all three solution methods, which are again referred to as Methods~1--3 as in Section~\ref{subsec:results-artificial}. The termination criteria based on running time limits and MIP gap were also adopted from Section~\ref{subsec:results-artificial}. The results, encompassing objective values and running time data, are presented in Table~\ref{tab:realWorldObjectives}.

\medskip

\begin{table}[ht]
    \centering
    \caption{Real-world instances: Running times and objective values} \label{tab:realWorldObjectives}
    \begin{adjustbox}{max width=\linewidth} % Adjust the table width
    \begin{tabular}{rrrrrrr}
        \toprule
        \multicolumn{1}{c}{\multirow{2}{*}{Instance~\#}}  & \multicolumn{2}{c}{\textbf{Method~1: MIP}} & \multicolumn{2}{c}{\textbf{Method~2: Sequ. sol. app.}} & \multicolumn{2}{c}{\textbf{Method~3: Heur. sol. app.}} \\
        & \multicolumn{1}{c}{Running time} & \multicolumn{1}{c}{Obj. value\tnote{1}} & \multicolumn{1}{c}{Running time} & \multicolumn{1}{c}{Obj. value\tnote{1}} & \multicolumn{1}{c}{Running time} & \multicolumn{1}{c}{Obj. value\tnote{1}} \\
        \midrule
            1 & 	10800 & 	100\% & 	4476 & 	103\% & 	10 & 	124\% \\ 
            2 & 	10800 & 	100\% & 	4346 & 	110\% & 	9 & 	129\% \\ 
            3 & 	10800 & 	100\% & 	3300 & 	106\% & 	9 & 	146\% \\ 
            4 & 	10800 & 	100\% & 	3142 & 	100\% & 	10 & 	134\% \\ 
            5 & 	10800 & 	100\% & 	4593 & 	107\% & 	8 & 	127\% \\ 
            6 & 	10800 & 	100\% & 	3397 & 	102\% & 	10 & 	127\% \\ 
            7 & 	10800 & 	100\% & 	1479 & 	106\% & 	11 & 	138\% \\ 
            8 & 	10800 & 	100\% & 	2792 & 	102\% & 	12 & 	128\% \\ 
            9 & 	10800 & 	100\% & 	1997 & 	112\% & 	10 & 	147\% \\ 
            10 & 	10800 & 	100\% & 	1894 & 	109\% & 	13 & 	140\% \\ 
            11 & 	10800 & 	100\% & 	1746 & 	109\% & 	10 & 	137\% \\ 
            12 & 	10800 & 	100\% & 	2358 & 	110\% & 	11 & 	164\% \\ 
            13 & 	10800 & 	100\% & 	294 & 	111\% & 	9 & 	148\% \\ 
            14 & 	10800 & 	100\% & 	113 & 	123\% & 	6 & 	142\% \\ 
            15 & 	10800 & 	100\% & 	656 & 	109\% & 	8 & 	121\% \\ 
            16 & 	10800 & 	100\% & 	4551 & 	104\% & 	8 & 	144\% \\ 
            17 & 	10800 & 	100\% & 	1920 & 	103\% & 	9 & 	122\% \\ 
            18 & 	10800 & 	100\% & 	403 & 	107\% & 	10 & 	125\% \\ 
            19 & 	10800 & 	100\% & 	1233 & 	109\% & 	11 & 	139\% \\ 
            20 & 	10800 & 	100\% & 	6006 & 	101\% & 	13 & 	161\% \\        
        \midrule
        \multicolumn{1}{c}{\textbf{Total}\tnote{2}} & 	10800 & 	100\% & 	2535 & 	107\% & 	10 & 	137\% \\ 
        \bottomrule
    \end{tabular}
    \end{adjustbox}
    \begin{tablenotes}
        \item \tnote{1} All objective values are provided as proportions (expressed as percentages) of those obtained using Method~1 (MIP).
        \item \tnote{2} Average across all instances
    \end{tablenotes}
\end{table}

The results in real-world instances parallel those in the artificial ones (Section~\ref{subsec:results-artificial}) concerning running time limits, where Method~1 consistently reaches the time limit in all cases. In Method~2, the running time limit is only exceeded in the PRA subproblem of instance~20. In stark contrast, Method~3 demonstrates remarkable efficiency, with an average running time of just 10~seconds, vastly outperforming Method~2, which has an average running time of 2535~seconds.

\medskip

Regarding objective values, Method~1 consistently outperforms Methods 2 and 3 across all instances, with average objective values of 107\% and 137\%, respectively. The real-world instances are less complex than the artificial 30~beds scenario, evidenced by an average MIP gap of 43\% (Table~\ref{tab:mipgap_realworld} in the appendix). Method~2 optimally solves all PRA subproblems except for instance~20, where a 100\%~MIP gap occurs because the lower bound equaled 0. In NPA subproblems, the MIP gap drops below 5\% during runtime, prompting the optimization to halt accordingly.

\medskip

In summary, the evaluation of three methods for healthcare management optimization reveals distinct trade-offs. Method~1, while suitable for small-scale settings such as those provided by our partner hospital, suffers from prohibitively high running time costs, often taking thousands of seconds to achieve results comparable to the heuristic, which only needs a few seconds. For medium-sized problems like the 30~beds scenario, Method~2 delivers superior results compared to Method~3 at the expense of significant running time. In contrast, Method~3 consistently stands out regarding running time efficiency across problem sizes and showcases excellent objective values, particularly in larger scenarios like the 60~beds scenario. Hence, Method~3 emerges as the preferred choice for optimizing healthcare management, striking a favourable balance between computational efficiency and solution quality, especially in more extensive and complex healthcare settings.

\section{Conclusion and outlook}\label{sec:conclusion}
Motivated by important interactions of PRA and NPA decisions in hospital wards, this paper explicitly considers both types of assignment decisions in one integrated optimization problem for the first time. We introduce the IPRNPA problem and provide a formal mathematical description as a mixed integer program. Since the PRA problem and the NPA problem are already \textsf{NP}-hard and very difficult to solve for realistic instance sizes, it is not surprising that the integrated problem is computationally challenging and cannot be solved to (near) optimality in reasonable time using the completely integrated MIP model. Therefore, we also present an efficient heuristic for the integrated problem that is able to compute high-quality solutions quickly on both artificially generated instances and real-world instances obtained from our partner hospital.
The heuristic solution approach highlighted its superiority in running time efficiency across various problem sizes, particularly excelling in larger and more complex scenarios. While the integrated MIP struggles with high running time costs the heuristic solution approach strikes a favorable balance between computational efficiency and solution quality. Its superiority in objective values is notably pronounced in medium and large-scale scenarios, making it the preferred choice for optimizing real-world healthcare settings.

\medskip
We also devise a parameterized instance generator for the problem. This generator is made freely available to other researchers to foster additional investigations on the IPRNPA problem, which we believe represents a challenging and at the same time practically relevant problem to be further investigated in the future. For instance, while our heuristic solution method allows easy adaptions to some dynamic versions of the problem, explicitly investigating different dynamic extensions with increasing degrees of data uncertainty (e.g., no-shows of patients or unexpected changes of patients' care requirements and/or LOS after admission) might represent a fruitful direction for future research. Moreover, while the nurse roster for the planning period is considered as an input of the problem in this paper, integrating rostering decisions into the problem formulation might represent an interesting extension.

% Future work could consider the dynamic admission of patients, e.g. emergency patients, who are transferred to the wards from the emergency department, or elective no-show patients, which could lead to changes in the room and nurse assignments. While the MIP only addresses the static version of the integrated problem in which all information about patients is known in advance, the heuristic can be easily adapted to dynamic settings. Doing do would make the heuristic even more relevant for an application in practice. The instance generator could be extend to design realistic dynamic instances to for an analysis of the differences and potential benefits between dynamic and static patient admission.

% Further research topics? Link to surgery scheduling maybe? Or to ICU management?
\section*{Statements and declarations}

\subsection*{Data availability}
The instance generator and the instances used for the computational experiments are available on GitHub, \url{https://github.com/TabeaBrandt/instance_generation_integrated_beds_and_staff_planning.git}.

\subsection*{Funding}
This research was funded by the Deutsche Forschungsgemeinschaft (DFG, German Research Foundation) -- Project number 443158418.

\smallskip

\noindent
Joe Viana's work is primarily supported by the Norwegian Research Council, Measure for Improved Availability of medicines and vaccines -- Project number 300867.

\subsection*{Conflicts of interest}
The authors have no conflicts of interest to declare that are relevant to the content of this article.
\bibliographystyle{sn-basic} 
\bibliography{bibfile}

\clearpage
\begin{appendices}
\section{Nurse rostering formulation}\label{sec:A1-nurse-rostering}
% \section{Nurse rostering formulation}

We use a simple binary integer programming formulation to generate the nurse rosters that are part of our random instances.
The formulation is based on the description presented in the first International Nurse Rostering Competition~(INRC) 2010~\cite{Haspeslagh2014}. As the focus of this work is not on nurse rostering, we use a simple but fast formulation instead of a perfectly detailed one. This formulation includes all constraints that are relevant concerning the use of a nurse roster as an input of the IPRNPA.
However, it is very easy to substitute the used nurse rostering formulation in our code.

\medskip

Similar to the INRC, we determine the roster for the planning period considering one ward. We use a subset of constraints of the INRC in order to compute a simple, yet still realistic nurse roster. These constraints include that the number of required nurses per shift must be met, not more than a given maximum allowed number of shifts can be assigned to any single nurse during the planning period, and that minimum rest times for nurses between shifts are respected. 

\medskip

\noindent
In addition to some of the notation and parameters introduced in Section~\ref{sec:MIP}, we use the following parameters and decision variables:
%$\skillnurses(s,l)$ which specifies for every shift $s\in \setS$ how many nurses of at least skill level $l \in \setL$ are needed.

\medskip

\noindent
\textbf{Parameters:}
\begin{itemize}[align=left,labelindent=0pt,labelwidth=2.7cm,labelsep*=1em,leftmargin=!]
\item[$\skillnurses(s,l)$\;] number of nurses with at least skill level $l \in \setL$ required during shift $s \in \setS$. %on day $d \in \setD$.
The sum of this number over all skill levels defines the minimum number of nurses needed per shift in total.
\item[$\maxshifts$\quad] maximum allowed number of shifts per nurse within the time horizon.
\end{itemize}

\medskip

\noindent
\textbf{Decision variables:}
\begin{itemize}[align=left,labelindent=0pt,labelwidth=2.7cm,labelsep*=1em,leftmargin=!]
\item[$\assign_{n,s}$\quad] binary variable indicating whether nurse~$n\in\setN$ is assigned to shift~$s\in\setS$
\end{itemize}

\noindent
The focus of our formulation is on the generation of a feasible roster to be used as an input for the IPRNPA. Therefore, we consider the minimization of the total number of assigned nurses as our objective in order to prevent unnecessary assignments that are not required in order to fulfill the considered constraints:

\medskip

\noindent
\textbf{Assignment objective}
\smallskip
\begin{enumerate}[resume,label=(\arabic*)]
\item Minimization of the number of assigned nurses: \label{objective:assignment}
\begin{align*}
\min \sum_{s\in\setS,n\in\setN} \assign_{n,s}
\end{align*}
\end{enumerate}

\noindent
The following constraints must be met by a feasible nurse roster:

\medskip

\noindent
\begin{enumerate}[resume,label=(\arabic*)]
\item Each nurse can work at most one shift per day: \label{constr_nurse:max1shift}
\begin{align*}
\assign_{n,s} + \assign_{n,s+1} + \assign_{n,s+2} \leq 1 & \quad\forall n\in\setN, s \in \setEarly
\end{align*}
\item For each skill level~$l \in \setL$, at least $\skillnurses(s,l)$ nurses with at least skill level~$l$ must be assigned during each shift~$s\in\setS$:
\label{constr_nurse:min-number-nurses-per-level}
\begin{align*}
\sum_{n\in\setN : \skilllevel(n) \geq l} \assign_{n,s} \geq \skillnurses(l,s) & \quad\forall l \in \setL, s\in\setS
\end{align*}
\item The minimum total number of nurses must be assigned during each shift~$s\in\setS$:
\label{constr_nurse:min-number-nurses}
\begin{align*}
\sum_{n\in\setN} \assign_{n,s} \geq \sum_{l \in \setL} \skillnurses(l,s) & \quad\forall s\in\setS
\end{align*}
\item No nurse~$n\in\setN$ can be assigned to more than~$\maxshifts$ many shifts during the planning period: \label{constr_nurse:max-number-shifts}
\begin{align*}
\sum_{s\in\setS} \assign_{n,s} \leq \maxshifts & \quad\forall n\in\setN
\end{align*}
\item On the day after a night shift, a nurse can only have another night shift (or the day off): \label{constr_nurse:after-night-shift}
\begin{align*}
%\assign_{n,s} + \assign_{n,s+2} \leq 1 & \quad\forall n\in\setN, s \in \setNight \\
%\assign_{n,s} + \assign_{n,s+1} \leq 1 & \quad\forall n\in\setN, s \in \setNight\\
\assign_{n,s} + \assign_{n,s+1} + \assign_{n,s+2} \leq 1 & \quad\forall n\in\setN, s \in \setNight
\end{align*}
\item On the day after a late shift, a nurse cannot have a morning shift: \label{constr_nurse:after-late-shift}
\begin{align*}
\assign_{n,s} + \assign_{n,s+2} \leq 1 & \quad\forall n\in\setN, s \in \setLate
\end{align*}
\end{enumerate}

\clearpage
\section{MIP Gap evaluation on artificial and real-world instances}\label{sec:B1-MIPgap}

\begin{table}[h]
    \centering
    \caption{Artificial instances: Average MIP gap} \label{tab:mipgap_artificial}
    \begin{adjustbox}{max width=\linewidth} % Adjust the table width
    \begin{tabular}{cccrrr}
        \toprule
        \multirow{2}{*}{Scenario} & \multirow{2}{*}{Variation} & \multirow{2}{*}{\shortstack{Planning\\horizon}} & \multirow{1}{*}{Method 1} & \multicolumn{2}{c}{Method 2} \\
        & & &\multicolumn{1}{c}{MIP gap\tnote{1}} & \multicolumn{2}{c}{MIP gap\tnote{1}} \\
        & & &\multicolumn{1}{c}{IPRNPA} & \multicolumn{1}{c}{PRA} & \multicolumn{1}{c}{NPA} \\
        \midrule
\multirow{6}{*}{30 beds} & \multirow{2}{*}{Var. 1 (double rooms)} & 2 weeks &	70.0\% & 	58.0\% & 	7.6\% \\ 
&   & 4 weeks &	80.7\% & 	63.4\% & 	10.9\% \\ 
			
& \multirow{2}{*}{Var. 2 (triple rooms)} & 2 weeks &	57.6\% & 	82.7\% & 	5.0\% \\ 
&   & 4 weeks &	80.8\% & 	79.7\% & 	6.8\% \\ 
			
& \multirow{2}{*}{Var. 3 (mixed rooms)} & 2 weeks &	65.3\% & 	55.0\% & 	6.1\% \\ 
&   & 4 weeks &	81.9\% & 	68.9\% & 	8.8\% \\ 
			
\midrule			
\multirow{6}{*}{60 beds} & \multirow{2}{*}{Var. 1 (double rooms)} & 2 weeks &	110.9\% & 	75.2\% & 	60.5\% \\ 
&   & 4 weeks &	111.1\% & 	89.3\% & 	62.3\% \\ 
			
& \multirow{2}{*}{Var. 2 (triple rooms)} & 2 weeks  &	98.0\% & 	94.6\% & 	51.1\% \\ 
&   & 4 weeks &	109.5\% & 	97.0\% & 	62.1\% \\ 
			
& \multirow{2}{*}{Var. 3 (mixed rooms)} & 2 weeks &	110.0\% & 	74.4\% & 	58.2\% \\ 
&   & 4 weeks &	110.3\% & 	90.5\% & 	61.4\% \\ 

        \bottomrule
    \end{tabular}
    \end{adjustbox}
    \begin{tablenotes}
        \item \tnote{1} The MIP gap is evaluated after either a running time limit of 10,800 seconds or 5,400 seconds for each subproblem or when the MIP gap reaches a value less than 5\%. The MIP gap is shown summarized as an average value.
    \end{tablenotes}
\end{table}
\end{appendices}

\begin{table}[ht]
    \centering
    \caption{Real-World instances: MIP gap} \label{tab:mipgap_realworld}
    \begin{adjustbox}{max width=\linewidth} % Adjust the table width
    \begin{tabular}{rrrr}
        \toprule
        \multicolumn{1}{c}{\multirow{3}{*}{Instance~\#}}  & \multirow{1}{*}{Method 1} & \multicolumn{2}{c}{Method 2} \\
       & \multicolumn{1}{c}{MIP gap\tnote{1}} & \multicolumn{2}{c}{MIP gap\tnote{1}} \\
        & \multicolumn{1}{c}{IPRNPA} & \multicolumn{1}{c}{PRA} & \multicolumn{1}{c}{NPA} \\
        \midrule
            1 & 	53\% & 	0\% & 	5\% \\ 
            2 & 	54\% & 	0\% & 	5\% \\ 
            3 & 	45\% & 	0\% & 	5\% \\ 
            4 & 	48\% & 	0\% & 	5\% \\ 
            5 & 	51\% & 	0\% & 	5\% \\ 
            6 & 	46\% & 	0\% & 	5\% \\ 
            7 & 	45\% & 	0\% & 	5\% \\ 
            8 & 	48\% & 	0\% & 	5\% \\ 
            9 & 	35\% & 	0\% & 	5\% \\ 
            10 & 	38\% & 	0\% & 	5\% \\ 
            11 & 	40\% & 	0\% & 	5\% \\ 
            12 & 	47\% & 	0\% & 	5\% \\ 
            13 & 	33\% & 	0\% & 	5\% \\ 
            14 & 	23\% & 	0\% & 	5\% \\ 
            15 & 	48\% & 	0\% & 	5\% \\ 
            16 & 	53\% & 	0\% & 	5\% \\ 
            17 & 	44\% & 	0\% & 	5\% \\ 
            18 & 	35\% & 	0\% & 	5\% \\ 
            19 & 	30\% & 	0\% & 	5\% \\ 
            20 & 	38\% & 	100\% & 	5\% \\       
        \midrule
        \multicolumn{1}{c}{\textbf{Total}\tnote{2}}  & 	43\% & 	5\% & 	5\% \\ 
        \bottomrule
    \end{tabular}
    \end{adjustbox}
    \begin{tablenotes}
        \item \tnote{1} The MIP gap is evaluated after either a running time limit of 10,800 seconds or 5,400 seconds for each subproblem or when the MIP gap reaches a value less than 5\%.
        \item \tnote{2} Average across all instances
    \end{tablenotes}
\end{table}

\end{document}